\documentclass[12pt]{article}
\usepackage{amssymb,amsmath,amsfonts}
\usepackage{relsize}
\usepackage[sort]{cite}
\usepackage{esvect}
\usepackage{txfonts}
\usepackage{graphicx}
\newcommand{\ver}{{\rm ver}}
\newcommand{\vo}{{\rm vol}}

\newcommand{\mes}{{\rm mes_n}}
\newcommand{\lin}{{\rm lin}}
\newcommand{\simp}{{\rm simp_n}}

\newcommand{\conv}{{\rm conv}}

\newcommand{\vol}{{\rm vol}}

\newcommand {\ip}[1]   {\langle{#1}\rangle}
\newcommand {\R}       {{\mathbb R}}
\newcommand {\RN}      {\R^n}
\newcommand {\smsk}    {\smallskip}
\newcommand {\msk}     {\medskip}
\newcommand {\bsk}     {\bigskip}

\newcommand {\VST}     {\vspace*{1mm}}

\newcommand {\SECT}[2] {\section*{\centerline{\normalsize
{\bf #1}}} \setcounter{section}{#2}
\setcounter{theorem}{0}\setcounter{equation}{0}}

\begin{document}
\parindent 1em
\parskip 0mm
\medskip
\centerline{\LARGE Geometric Estimates in Linear Interpolation}
\vspace*{3mm}
\centerline{\LARGE on a Cube and a Ball}\vspace*{10mm}
\centerline{\large Mikhail Nevskii}\vspace*{5 mm}
\centerline {\it Department of Mathematics,  P.\,G.~Demidov Yaroslavl State University,}
\centerline {\it Sovetskaya str., 14, Yaroslavl, 150003, Russia}\vspace*{5 mm}
\centerline{\large February 18, 2024}

\vspace*{15mm}
\renewcommand{\thefootnote}{ }
\footnotetext[1]{\hspace{-6mm}
{\it E-mail address:} mnevsk55@yandex.ru}
\hrule\msk

\bsk
\par\noindent {\bf Abstract}
\bsk
\par The paper contains a survey of the results obtained by the author in recent years.
These results concern the application in multivariate polynomial interpolation of some geometric constructions and methods. In particular, we give estimates of the projector's norms through the characteristics of sets associated with homothety.
The known exact values and nowaday best estimates of these norms are~given.
Also we formulate some open problems.
\bsk
{\small
\smsk
}
\msk
\hrule

\renewcommand{\contentsname}{ }
\tableofcontents
\addtocontents{toc}{{\bf Contents}\vspace*{5mm}\par}
\bsk


\SECT{Introduction}{1}\label{intro}
\addtocontents{toc}{Introduction\hfill \thepage\par\VST}
\indent

The results related to minimization the norms of projection operators (briefly, projectors) are important in Approximation Theory thanks to the inequality
$$
E_m(f;K)_{C(K)}\le \|f-Pr(f)\|_{C(K)}\le (1+\|Pr\|)E_m(f;K)_{C(K)}.
$$
Here $E_m(f;K)_{C(K)}$ is the best approximation in the uniform norm of a continuous function $f$ defined on a closed set $K\subset\RN$ by polynomials of degree at most $m$ (the space $\Pi_m(\RN)$), $Pr:C(K)\to \Pi_m(\RN)$ is a bounded projector from $C(K)$ onto $\Pi_m(\RN)$, and
$\|Pr\|$ is its operator norm.
\par The paper contains  a survey of the results obtained by the author in recent years. They are related to geometric aspects of multidimensional polynomial interpolation. The most deep results concern linear interpolation.
In particular, we give estimates of the projector's norms through the characteristics of~sets associated with homothety.
The known exact values and nowaday best estimates of~these norms are given. The approach related to linear interpolation can also be used for interpolation by~polynomials from wider spaces. Our main focus is on a cube and a ball in ${\mathbb  R}^n$. A number of~statements we give with proofs, otherwise references to the literature are given. Also we formulate some open questions.

Numerical characteristics connecting
simplices and subsets of ${\mathbb R}^n$
have applications for obtaining various estimates
in polynomial interpolation.
This approach and the corresponding analytic
methods
in details
were described in the author's works
\cite{nevskii_mais_2003_10_1}--\cite{nevskii_fpm_2013} including the monograph \cite{nevskii_monograph}.
 Lately these questions have been managed to study
also by computer methods (see \cite{nev_ukh_mais_2016_23_5}, \cite{nev_ukh_mais_2017_24_1},
\cite{nev_ukh_mais_2018_25_1}, \cite{nev_ukh_beitrage_2018}, \cite{nev_ukh_mais_2018_25_3},
\cite{nev_ukh_saratov_2018},
\cite{nev_ukh_mais_2019_26_2}, \cite{nev_ukh_matzam SVFU_2019},
\cite{nev_ukh_dataset_2019}, \cite{nev_ukh_mais_2023_30_}, \cite{ukhalov_udovenko_2020},
and also the survey paper \cite{nevskii_ukh_sb_2021}).
Computer calculations were carried out by Alexey Ukhalov, whose skill and persistence produced many impressive results. Various computer math systems, such as Wolfram Mathematica (see, e.\,g., \cite{wolfram_2016}),
 were used significantly along the way. Some of the numerical results were obtained by our students.

 \bigskip
 In the early 1970s, at Yaroslavl' State University (Yaroslavl', Russia) there was formed a~school on Function Theory, Approximation Theory, Linear Operator Interpolation Theory, and related areas.  The founder and head of this school for about 20 years was Professor Yuri Abramovich Brudnyi. It~was Prof. Brudnyi who inspired the author to take up Approximation Theory. Under his supervision the author performed  his dissertation research related to piecewise-polynomial approximation in~Orlicz classes.

 Throughout their lives, all of
Yu.\,A. Brudnyi's students have felt great gratitude, love and respect for their teacher. This text is dedicated to Yuri Abramovich Brudnyi in~connection with the  upcoming  90th anniversary of his birth.

\SECT{1. Notation and preliminaries}{1}
\label{main_defs_and_denots}
\addtocontents{toc}{1. Notation and preliminaries\hfill \thepage\par\VST}


\indent\par Let us fix some notation. Throughout the paper  $n\in{\mathbb N}$ is a positive integer.
By $e_1$, $\ldots$, $e_n$ we denote the standard basis
in ${\mathbb R}^n$. Given $x=(x_1,...x_n)\in\RN$, by $\|x\|$ we denote its standard Euclidean norm
$$
\|x\|=\sqrt{\langle x,x\rangle}=\left(\sum\limits_{i=1}^n x_i^2\right)^{1/2}.
$$
Hereafter, for $x=(x_1,...x_n), y=(y_1,...y_n)\in\RN$ by
$
\ip{x,y}=x_1y_1+...+x_ny_n
$
we denote the standard inner product in $\RN$.

\par Given $x^{(0)}\in\RN$ and $R>0$, we let $B(x^{(0)};R)$ denote the Euclidean ball with center $x^{(0)}$ and radius $R$:
$$
B\left(x^{(0)};R\right)=\{x\in{\mathbb R}^n: \|x-x^{(0)}\|\leq R \}.
$$

We also set
$$
B_n=B(0;1),~~~Q_n=[0,1]^n~~~\text{and}~~~
Q_n^\prime:=[-1,1]^n.
$$

The notation $L(n)\asymp M(n)$ means that there exist absolute constants $c_1,c_2>0$ such that 
\linebreak $c_1M(n)\leq L(n)\leq c_2 M(n)$.

Let $K$ be {\it a convex body in ${\mathbb R}^n$}, i.e., a compact convex subset of ${\mathbb R}^n$
with nonempty interior. By $c(K)$ we denote the center of gravity of $K$.
The symbol $\sigma K$  denotes a homothetic copy of $K$ with the center of homothety in $c(K)$
and the ratio of homothety  $\sigma.$

We let $\vol(K)$ denote the volume of $K$.
If $K$ is a convex polytope, then by $\ver(K)$ we denote the set of all vertices of $K$. By {\it a translate} we mean the result of a parallel shift.

We say that an $n$-dimensional simplex $S$ is {\it circumscribed around a convex body $K$} if $K\subset S$ and~each $(n-1)$-dimensional face of  $S$ contains a point of $K$. A convex polytope is {\it inscribed into $K$} if every vertex of this polytope belongs to the boundary of
$K$.

\smallskip
{\bf Definition 1.1.} {\it Let $i\in\{1,...,n\}$ and let $d_i(K)$ be the maximal length of a segment contained in $K$ and parallel to the
$x_i$-axis. We refer to $d_i(K)$  as {\it the $i$th
axial diameter of $K$.
}
}

\smallskip
The notion of {\it axial diameter of a convex body} was introduced by P.\,Scott \cite{scott_1985},\cite{scott_1989}.

\smallskip
{\bf Definition 1.2.} {\it  Given convex bodies $K_1$, $K_2$, by $\xi(K_1;K_2)$ we denote the minimal $\sigma\geq 1$ with the~property
$K_1\subset \sigma K_2$.  We call $\xi(K_1,K_2)$ {\it the absorption index of $K_1$ by $K_2$}.

By $\alpha(K_1,K_2)$ we denote the minimal $\sigma>0$ such that {\it $K_1$ is a subset of a translate of $\sigma K_2$}.
}

\smallskip
Note that the equality
$\xi(K_1;K_2)$ $=1$ is equivalent to the inclusion  $K_1\subset K_2$. Clearly, $\alpha(K_1,K_2)$ $\leq $ $\xi(K_1,K_2)$.

\smallskip
{\bf Definition 1.3.} {\it By $\xi_n(K)$ we denote the minimal absorption index of a convex body $K$ by an inner nondegenerate simplex.
In other words,
$$
\xi_n(K)=\min \{ \xi(K;S): \,
S  \mbox{ is an $n$-dimensional simplex,} \,
S\subset K, \, \vo(S)\ne 0\}.
$$}

\smallskip
By $C(K)$ we denote the space of all continuous functions
$f:K\to{\mathbb R}$ with the uniform
norm
$$
\|f\|_{C(K)}=\max\limits_{x\in K}|f(x)|.
$$
We let $\Pi_1\left({\mathbb R}^n\right)$ denote the space of polynomials in $n$ variables
of degree at most $1$.

Let $S$ be a nondegenerate simplex in ${\mathbb R}^n$ with vertices $x^{(j)}=\left(x_1^{(j)},\ldots,x_n^{(j)}\right),$
$1\leq j\leq n+1.$ We define {\it the vertex matrix} of this simplex by
$$
{\bf A}=
\left( \begin{array}{cccc}
x_1^{(1)}&\ldots&x_n^{(1)}&1\\
x_1^{(2)}&\ldots&x_n^{(2)}&1\\
\vdots&\vdots&\vdots&\vdots\\
x_1^{(n+1)}&\ldots&x_n^{(n+1)}&1\\
\end{array}
\right).
$$

Clearly, matrix ${\bf A}$ is nondegenerate and
\begin{equation}\label{vol_S_det_A_eq}
\vol(S)=\frac{|\det({\bf A})|}{n!}.
\end{equation}
Let
\begin{equation}\label{ACV}
{\bf A}^{-1}=\left(l_{ij}\right)_{i,j=1}^{n+1}.
\end{equation}

\smallskip
{\bf Definition 1.4.} {\it Linear  polynomials
$$
\lambda_j(x)=
l_{1j}x_1+\ldots+
l_{nj}x_n+l_{n+1,j},~~~~j=1,...,n+1,
$$
are called {\it the~basic Lagrange polynomials} corresponding to $S$.
}

\smallskip
These polynomials have the following property:
$$
\lambda_j\left(x^{(k)}\right)=\delta_j^k~~~\text{for all}~~~j,k=1,...,n+1.
$$
Here $\delta_j^k$ is the Kronecker delta.
For an arbitrary $x\in{\mathbb R}^n$, we have
$$
x=\sum_{j=1}^{n+1} \lambda_j(x)x^{(j)}, \quad \sum_{j=1}^{n+1} \lambda_j(x)=1.
$$
Thus, $\lambda_j(x)$ are {\it the barycentric coordinates of $x$} with respect to the simplex $S$.
In turn, equations $\lambda_j(x)=0$, $j=1,...,n+1$, define the $(n-1)$-dimensional hyperplanes containing the faces of $S.$  Therefore,
$$
S=\left \{ x\in {\mathbb R}^n: \, \lambda_j(x) \geq 0, \, j=1,\ldots,n+1
\right\}.
$$

Also let us note that for every $j=1,...,n+1$ we have
\begin{equation}\label{Lagr_pol_thru_dets}
\lambda_j(x)=\frac{\Delta_j(x)}{\Delta}.
\end{equation}
Here $\Delta=\det({\bf A})$ and $\Delta_j(x)$ is the determinant that appears from $\Delta$ after replacing the $j$th row by~the~row $(x_1\  \ldots\  x_n \ 1).$
For more information on $\lambda_j$, see ~\cite{nevskii_monograph}, \cite{nev_ukh_posobie_2020}.

\smallskip
{\bf Definition 1.5.} {\it Let $x^{(j)}\in K$, $1\leq j\leq n+1,$ be  the vertices of a nondegenerate simplex $S$.
The interpolation projector $P:C(K)\to \Pi_1({\mathbb R}^n)$  with nodes $x^{(j)}$ is determined by the equalities $Pf\left(x^{(j)}\right)=f_j=f\left(x^{(j)}\right)$ , $1\leq j\leq n+1.$
We say that an interpolation projector $P:C(K)\to \Pi_1({\mathbb R}^n)$ and a nondegenerate
simplex $S\subset K$ {\it correspond to each other} if
the  nodes of $P$ coincide with  the~vertices of $S$. We use notation $P_S$ and $S_P$.
}

\smallskip
For an interpolation projector $P=P_S$, the analogue of Lagrange interpolation formula holds:
\begin{equation}\label{interp_Lagrange_formula}
Pf(x)=\sum\limits_{j=1}^{n+1}
f\left(x^{(j)}\right)\lambda_j(x).
\end{equation}
Here $\lambda_j$ are the basic Lagrange polynomials of the simplex $S_P$ (see Definition 1.4).

We let $\|P\|_K$ denote the norm of $P$ as an operator from $C(K)$ into $C(K)$.  Thanks to (\ref{interp_Lagrange_formula}),
$$
\|P\|_K=\sup_{\|f\|_{C(K)}=1} \|Pf\|_{C(K)}=
\sup_{-1\leq f_j\leq 1} \max_{x\in K}\left| \sum_{j=1}^{n+1} f_j\lambda_j(x)\right|=$$
$$=
\max_{x\in K}\sup_{-1\leq f_j\leq 1}\left| \sum_{j=1}^{n+1} f_j\lambda_j(x)\right|=
\max_{x\in K}\sup_{-1\leq f_j\leq 1}\sum_{j=1}^{n+1} f_j\lambda_j(x).$$
Because $\sum f_j\lambda_j(x)$ is linear in $x$ and
$f_1,\ldots,f_{n+1}$, we have
\begin{equation}\label{norm_P_intro_cepochka}
\|P\|_K=
\max_{x\in K}\max_{f_j=\pm 1} \sum_{j=1}^{n+1}
f_j\lambda_j(x)
=\max_{x\in K}\sum_{j=1}^{n+1}
|\lambda_j(x)|.
\end{equation}
If $K$ is a convex polytope in ${\mathbb R}^n$
(e.\,g., $K$ is a cube), a simpler equality
\begin{equation}\label{norm_P_cube_formula}
\|P\|_K= \max_{x\in\ver(K)}\sum_{j=1}^{n+1}
|\lambda_j(x)|
\end{equation}
holds.

\smallskip
{\bf Definition 1.6.} {\it We let $\theta_n(K)$ denote the minimal value of $\|P_S\|_K$ where $S$ runs over all nondegenerate simplices with vertices in $K$. An interpolation projector $P:C(K)\to\Pi_1\left({\mathbb R}^n\right)$ is called minimal if $\|P\|_K=\theta_n(K)$. }

\smallskip
{\bf Definition 1.7.} {\it  Let $\mu$ be a positive integer, $1\leq \mu \leq n$. A point $x\in K$ is called
{\it the $\mu$-point of~$K$ with respect to a simplex $S\subset K,$} if the corresponding projector
$P:C(K)\to\Pi_1\left({\mathbb R}^n\right)$
satisfies $\|P\|_K=\sum | \lambda_j(x)|$
and, in addition,
among the numbers $\lambda_j(x)$ there are exactly $\mu$ negatives.

Suppose $K$ is a cube in $\RN$.  A  vertex $x$ of the cube is called  the $\mu$-vertex with respect to $S,$
if $x$ is~a~$\mu$-point of $K$ with respect to this  simplex.
}

\smallskip
It was shown in \cite{nevskii_mais_2008_15_3} that for any interpolation projector
$P:C(K)\to\Pi_1\left({\mathbb R}^n\right)$ and the corresponding simplex  $S$ we have
\begin{equation}\label{nev_ksi_P_ineq}
\frac{n+1}{2n}\Bigl( \|P\|_K-1\Bigr)+1\leq
\xi(K;S) \leq
\frac{n+1}{2}\Bigl( \|P\|_K-1\Bigr)+1.
\end{equation}

If there exists an $1$-point in $K$ with respect to $S$, then
the right-hand inequality in (\ref{nev_ksi_P_ineq}) becomes an equality.
This statement was proved in \cite{nevskii_mais_2008_15_3}
in the equivalent form; the notion of $1$-vertex was introduced later
in \cite{nevskii_mais_2009_16_1}. Furthermore, it was shown in \cite{nevskii_mais_2009_16_1} that if $K=Q$ where $Q$ is a {\it cube} in $\RN$ then the existence of an $1$-vertex of $Q$ with respect to a simplex $S\subset Q$ is equivalent
to the right-hand equality in \eqref{nev_ksi_P_ineq}.
 If
there exists a $\mu$-vertex of
$Q$ with respect to $S,$ then
$$
\frac{n+1}{2\mu}\left( \|P\|_{Q}-1\right)+1\leq \xi(Q;S).
$$

Thanks to \eqref{nev_ksi_P_ineq},
\begin{equation}\label{nev_ksi_n_K_theta_n_K_ineq}
\frac{n+1}{2n}\Bigl( \theta_n(K)-1\Bigr)+1\leq
\xi_n(K) \leq
\frac{n+1}{2}\Bigl( \theta_n(K)-1\Bigr)+1.
\end{equation}
Obviously, if a projector $P$ satisfies the equality
\begin{equation}\label{ksi_n_K_norm_P_eq}
\xi_n(K) =
\frac{n+1}{2}\Bigl( \|P\|_K-1\Bigr)+1,
\end{equation}
then $P$ is minimal and the right-hand relation in \eqref{nev_ksi_n_K_theta_n_K_ineq} becomes an equality.


It turns out that calculation of  $\xi(K;S)$ and $\alpha(K;S)$ is of particular interest to us.
In \cite{nevskii_dcg_2011} we show that
\begin{equation}\label{xi_K_S_formula}
\xi(K;S)=(n+1)\max_{1\leq k\leq n+1}
\max_{x\in K}(-\lambda_k(x))+1 \quad (K\not\subset S),
\end{equation}
\begin{equation}\label{alpha_K_S_general_formula}
\alpha(K;S)=\sum_{j=1}^{n+1} \max_{x\in K} (-\lambda_j(x))+1
\end{equation}
(see also \cite{nevskii_monograph}). Furthermore, we prove that the equality $\xi(K;S)=\alpha(K;S)$  holds true if and only if the simplex $\xi(S)S$ is circumscribed around $K$. This  is also equivalent to the relation
\begin{equation}\label{xi_S_S_circ_condition_around_K}
 \max_{x\in K} \left(-\lambda_1(x)\right)=
\ldots=
\max_{x\in K} \left(-\lambda_{n+1}(x)\right).
\end{equation}

If $K$ is a convex polytope, then the maxima on $K$ in 
\eqref{xi_K_S_formula}--\eqref{xi_S_S_circ_condition_around_K} can also be taken over $x\in \ver(K)$. Note that for $K=Q_n$
formula \eqref{xi_K_S_formula} is proved in \cite{nevskii_matzam_2010}.

Occasionally, we will consider the case when $n+1$ is~an~Hadamard number, i.e., there exists an~Hadamard matrix of~order $n+1$.

\smallskip
{\bf Definition 1.8.} {\it An  Hadamard matrix of order $m$ is a square binary matrix $\bf H$ with entries either $1$ or $-1$ satisfying the  equality
$${\bf H}^{-1}=\frac{1}{m}\, {\bf H}^{T}.$$
An integer $m$, for which an Hadamard matrix of order $m$ exists, is called an Hadamard number.}

\smallskip
Thus, the rows of ${\bf H}$ are pairwise orthogonal with respect to the standard scalar product on ${\mathbb R}^m.$

Let us recall some basic facts relating to the theory of Hadamard matrices. The order \linebreak of an Hadamard matrix is $1$ or $2$ or some multiple of $4$ (see \cite{hall_1970}).  But it is still unknown whether an Hadamard matrix exists for every order of the form $m=4k$. This  is one of the longest lasting open problems in Mathe\-matics called {\it the Hadamard matrix conjecture}.
The orders below $1500$ for which Hadamard matrices are not known are $668, 716, 892, 956, 1132, 1244, 1388$, and $1436$ (see, e.g., \cite{horadam_2007}, \cite{manjhi_2022}).

Two Hadamard matrices are called equivalent if one can be obtained  from another by performing a finite sequence of the following operations: multiplication of~some rows or columns by $-1$, or~permutation of rows or columns. Up to equi\-valence, there exists a unique Hadamard matrix of orders $1, 2, 4, 8$, and $12$. There are $5$ equivalence classes of Hadamard matrices  of order $16, 3$ of order $20$,  $60$ of order $24$, and $487$ of order $28$.
For orders $32$, $36$, and $40$, the number of equivalence classes is much greater.
For $n=32$ there are at least $3,578,006$ equivalence classes; for $n=36$, at least $4,745,357$ (see \cite{horadam_2007}).

It is known that $n+1$ is an Hadamard number if and only if  it is possible to inscribe an $n$-dimensional regular simplex into an $n$-dimensional cube in such a way that the vertices of the simplex coincide with vertices of the cube. See \cite{hudelson_1996}.

Let us prove this statement. Let $Q_n^\prime=[-1,1]^n$. Then there exists a simple correspondence between Hadamard matrices of order $n+1$ with the last column consisting of $1$'s and $n$-dimensional regular simplices whose  vertices are at the vertices of $Q_n^\prime$. Since later on we will need this correspondence, we will dwell on it in more detail.

If $S$ is a regular simplex inscribed into  $Q_n^\prime$,
then its vertex matrix  ${\bf A}$  is an Hadamard matrix of order $n+1$. Indeed,
let  $x^{(1)}, \ldots,  x^{(n+1)}$ be the vertices of  $S$.
Since $x^{(j)}$, $j=1,...,n+1,$ coincide with vertices of the cube, the entries of $\bf A$ are  $\pm 1$ and~its last column consists of $1$'s.
Let us denote by $h^{(j)}$  the rows of  ${\bf A}$ and make sure that these vectors are pairwise orthogonal in ${\mathbb R}^{n+1}$.
Equalities $||x^{(j)}||= \sqrt{n}$ mean that simplex $S$ is inscribed into an $n$-dimensional ball of radius $\sqrt{n}$.
The edge-length $d$ of~a~regular simplex and radius $R$ of the circumscribed ball satisfy the following equality:
\begin{equation}\label{d_R_eq}
d=R\sqrt{2}\sqrt{\frac{n+1}{n}}
\end{equation}
 (see, e.\,g.,~\cite{nevskii_monograph}).
If $R=\sqrt{n},$ then $d=\sqrt{2(n+1)}.$  Hence,
$$2(n+1)=||x^{(j)}-x^{(k)}||^2=2||x^{(j)}||^2
- 2\langle x^{(j)},x^{(k)}\rangle=2n-2\langle x^{(j)},x^{(k)}\rangle.$$
 We see that  $\langle x^{(j)},x^{(k)}\rangle=-1.$
Since the rows of $\bf S$ in the first $n$ columns contain coordinates of vertices of the simplex and the last element of any row is $1$, we have
$$\langle h^{(j)},h^{(k)} \rangle=\langle x^{(j)},x^{(k)}\rangle+1=0.$$
Thus,  $h^{(j)}$  are pairwise orthogonal in ${\mathbb R}^{n+1}$ so that $\bf A$ is an Hadamard matrix of order $n+1$. Hence,
\begin{equation}\label{bfS_is_Hadamard}
{\bf A }^{-1}=\frac{1}{n+1}\, {\bf A}^{T}.
\end{equation}

Conversely, let ${\bf H}$ be an Hadamard matrix of order $n+1$ with the last column consisting of  $1$'s. Consider the simplex   $S$ whose vertex matrix is  $\bf H$.
In other words, the~vertices of  $S$ are given by~the~rows of   ${\bf H}$  (excepting the last component).  Let us see that the simplex $S$ is inscribed into the cube $Q_n^\prime,$ and~$\ver(S)\subset \ver(Q_n^\prime)$. Indeed, if $h^{(j)}$ are the rows of ${\bf H}$ and $x^{(j)}$ are the vertices of $S$, then
$\langle h^{(j)},h^{(k)}\rangle =0,$  $\langle x^{(j)},x^{(k)}\rangle =-1,$  $\|x^{(j)}\|^2=n.$
Hence, $||x^{(j)}-x^{(k)}||^2=2(n+1)$ proving that $S$ is a regular simplex with the edge-length  $\sqrt{2(n+1)}$.

Notice that any Hadamard matrix of order  $n+1$ is equivalent to the vertex matrix of some   $n$-dimensional regular simplex inscribed in $Q_n^\prime$. This vertex matrix
can be obtained from the given Hadamard  matrix after multiplication of some rows by~$-1$.

\smallskip
{\bf Definition 1.9.} {\it By $h_n$ we denote the maximum value of a determinant of order $n$ with entries $0$ or~$1$. Let $\nu_n$ be the maximum volume of an $n$-dimensional simplex
 contained in $Q_n$.}

 \smallskip
 The numbers $h_n$ and $\nu_n$ satisfy  the equality $h_n=n!\nu_n$ (see \cite{hudelson_1996}).
For any $n$, there exists  in $Q_n$ a~maximum volume simplex with some vertex coinciding with a vertex of the cube. 
For $n>1$, the following inequalities hold
\begin{equation}\label{adamar_clements_lindstrem1}\frac{1}{2}\left( 1-\frac{\log(4/3)}{\log n}\right ) n \log n
< \log(2^{n-1}h_{n-1})\leq \frac{1}{2} \, n\log n.
\end{equation}
The right-hand inequality in
(\ref{adamar_clements_lindstrem1}) was proved by Hadamard
\cite{hadamard_1893} and the left-hand one by Clements and  Lindstr\"om \cite{clements_lindstrem_1965}.
Consequently, for all $n\in {\mathbb N}$
\begin{equation}\label{adamar_clements_lindstrem2}
\left(\frac{3}{4}\right)^{(n+1)/2}\,\frac{\left(n+1\right)^{(n+1)/2}}{2^n}<
h_n\leq\frac{\left(n+1\right)^{(n+1)/2}}{2^n},
\end{equation}
\begin{equation}\label{adamar_clements_lindstrem}
\left(\frac{3}{4}\right)^{(n+1)/2}\,\frac{\left(n+1\right)^{(n+1)/2}}{2^nn!}<
\nu_n\leq\frac{\left(n+1\right)^{(n+1)/2}}{2^nn!}.
\end{equation}

The right-hand equality in each relation
holds if and only if $n+1$ is an Hadamard number \cite{hudelson_1996}.
In some cases the right-hand inequality in
(\ref{adamar_clements_lindstrem})
has been improved.
For instance, if $n$ is even, then
\begin{equation}\label{nu_n_n_even}
\nu_n\leq \frac{n^{n/2}\sqrt{2n+1}}{2^nn!}.
\end{equation}

If $n>1$ and $n\equiv 1({\rm mod}~4),$ then
\begin{equation}\label{nu_n_n_odd}
\nu_n\leq \frac{(n-1)^{(n-1)/2}}{2^{n-1}(n-1)!}
\end{equation}
(see \cite{hudelson_1996}).
For many $n$, the  values of $\nu_n$ and $h_n$ are known exactly.  The first $12$ numbers $\nu_n$ are
$$
\nu_1=1, \quad \nu_2= \frac{1}{2}, \quad \nu_3= \frac{1}{3},
\quad \nu_4=\frac{1}{8}, \quad \nu_5=\frac{1}{24}, \quad
\nu_6=\frac{1}{80}, \quad \nu_7=\frac{2}{315},$$\
$$\nu_8=\frac{1}{720}, \quad \nu_9=\frac{1}{2520}, \quad
\nu_{10}=\frac{1}{11340}, \quad \nu_{11}=\frac{9}{246400},
\quad \nu_{12}=\frac{3}{394240}.
$$

\medskip

\SECT{2. The values of $\alpha(Q_n;S)$ and $\xi(Q_n;S)$}{2}
\label{est_alpha_ksi_cube}
\addtocontents{toc}{2. The values of $\alpha(Q_n;S)$ and $\xi(Q_n;S)$ \hfill \thepage\par\VST}

\indent\par Recall that $d_i(K)$ is the maximal length of a segment contained in $K$ and parallel to the $x_i$-axis (see Definition 1.1).
 In \cite{nevskii_dcg_2011} we have proved the following statement.

\smallskip
{\bf Theorem 2.1.}  {\it Any convex body $K\subset\RN$ contains a translate of the cube $\sigma Q_n$ where
$$
\sigma=\left(\sum_{i=1}^n
\frac{1}{d_i(K)}
\right)^{-1}.
$$
}

This result immediately implies the  inequality
\begin{equation}\label{alpha_d_i_K_ineq}
\alpha(Q_n;K)\leq\sum_{i=1}^n\frac{1}{d_i(K)}.
\end{equation}
See  Definition 1.2
for $\alpha(Q_n;K)$.
It was also shown in \cite{nevskii_dcg_2011} that if $K=S$ where  $S$ is a nondegenerate simplex, then inequality \eqref{alpha_d_i_K_ineq} becomes an equality.

\smallskip
{\bf Theorem 2.2.} {\it For an arbitrary $n$-dimensional nondegenerate simplex $S$,
\begin{equation}\label{alpha_d_i_S_eq}
\alpha(Q_n;S)=\sum_{i=1}^n\frac{1}{d_i(S)}.
\end{equation}
}

\smallskip
Earlier, in \cite{nevskii_matzam_2010}, we have proved that
\begin{equation}\label{d_i_l_ij_formula}
\frac{1}{d_i(S)}=\frac{1}{2}\sum_{j=1}^{n+1} |l_{ij}|.
\end{equation}
Recall that $l_{ij}$, $i,j=1,...,n+1$, are the entries of the matrix ${\bf A}^{-1}$, see \eqref{ACV}.
Also we have shown that the only segment of length $d_i(S)$ parallel to the $i$th axis is located in $S$ in such a way that every $(n-1)$-dimensional face of the simplex contains at least one of its endpoints.
Formulas for the endpoints of this segment are also given in \cite{nevskii_matzam_2010}. For a generalization of these results to the case of maximum segments of {\it arbitrary directions} in a simplex, we refer the reader to \cite{nevskii_fpm_2013}.

\smallskip
{\bf Corollary 2.1.} {\it We have
\begin{equation}\label{alpha_qs_formula}
\alpha(Q_n;S)=\frac{1}{2}\sum_{i=1}^n\sum_{j=1}^{n+1} |l_{ij}|.
\end{equation}
}

\smallskip
This is immediate from \eqref{alpha_d_i_S_eq} and \eqref{d_i_l_ij_formula}.
Recall that $Q_n^\prime=[-1,1]^n.$

\smallskip
{\bf Corollary 2.2.} {\it
\begin{equation}\label{alpha_for_Q_prime}
\alpha(Q_n^\prime;S)=\sum_{i=1}^n\sum_{j=1}^{n+1} |l_{ij}|.
\end{equation}}

\smallskip
{\it Proof.} Obviously, $\alpha(K_1;K_2)$ is invariant with respect to shifts, and
$\alpha(\tau K_1;K_2)=\tau\alpha(K_1;K_2)$ for every $\tau>0$.  Since
$Q_n^\prime=[-1,1]^n$
is a translate of $2Q_n$, then \eqref{alpha_for_Q_prime} is the direct consequence of (\ref{alpha_qs_formula}).
\hfill$\Box$

\smallskip
We note that formula \eqref{alpha_d_i_S_eq} has many  useful corollaries which we present in works
\cite{nevskii_dcg_2011} and \cite{nevskii_monograph}. Let us give some examples.

\smallskip
{\bf Corollary 2.3.} {\it If  $S\subset Q_n$, then $\xi(Q_n;S)\geq n.$ The equality $\xi(Q_n;S)= n$
implies $\xi(Q_n;S)=\alpha(Q_n;S)$ and $d_i(S)=1$ for every $i=1,...,n$. }

\smallskip
{\it Proof.} Since $S\subset Q_n,$ we have $d_i(S)\leq 1$. Therefore,
\begin{equation}\label{xi_S_alpha_S_geq_n}
\xi(Q_n;S)\geq \alpha(Q_n;S)=\sum_{i=1}^n\frac{1}{d_i(S)}\geq n.
\end{equation}
If, in addition, $\xi(Q_n;S)= n$, then \eqref{xi_S_alpha_S_geq_n} gives  $\xi(Q_n;S)=\alpha(Q_n;S)=n$ and $d_i(S)=1$.
\hfill$\Box$

\smallskip
Equality \eqref{alpha_d_i_S_eq} yields the following property first  proven in \cite{nevskii_matzam_2010}.

\smallskip
{\bf Corollary 2.4.}
{\it If an $n$-dimensional simplex $S$ contains the cube $Q_n$, then for some $i\in\{1,...,n\}$ this~simplex contains a segment of length $n$ parallel to the $i$th axis.}

{\it Proof.} Because $Q_n\subset S$, we have $\alpha(Q_n;S)\leq 1.$ From \eqref{alpha_d_i_S_eq}, it follows
$$
\sum_{i=1}^n \frac{1}{d_i(S)}\leq 1
$$
proving the existence of $i$ such that $d_i(S)\geq n$.
\hfill$\Box$

\smallskip
{\bf Corollary 2.5.}
{\it Let   $S$ be a nondegenerate simplex and let
$D,$ $D^*$~be parallelotopes in ${\mathbb R}^n.$ Suppose
$D^*$ is a homothetic copy of $D$
with ratio $\sigma>1.$
If $D\subset S \subset D^*,$
then $\sigma\geq n.$}

\smallskip
{\it Proof.} It is sufficient to consider the case $D=Q_n.$ Then $D^*$ is a translate of the cube
$\sigma Q_n$. Since $S$ is contained in this translate, we have $d_i(S)\leq \sigma.$ Hence,
$1/\sigma \leq 1/d_i(S)$ for all $i$. Summing up over $i=1,\ldots,n$, we get
$$\frac{n}{\sigma}\leq\sum_{i=1}^n \frac{1}{d_i(S)}=\alpha(Q_n;S).$$
The inclusion $S\subset Q_n$ means that $\alpha(Q_n;S)\leq 1$ (for $\alpha(Q_n;S)$, see Definition 1.2).  Consequently, $n/\sigma \leq 1,$ i.\,e., $\sigma\geq n.$
\hfill$\Box$

\smallskip
Finally, let us note the  result  obtained earlier by M.\,Lassak \cite{lassak_dcg_1999}.

 \smallskip
{\bf Corollary 2.6.} {\it If a simplex $S\subset Q_n$ has the maximum volume, then $d_i(S)=1$ for every $i=1,...,n$.}

\smallskip
{\it Proof.} If a simplex $S\subset Q_n$ has the maximum volume, then $S\subset -nS$. In fact, if this is not the~case, then some vertex of $S$ can be moved in $Q_n$ in such a way that its distance from the opposite face of~the~simplex increases. In this case the volume of $S$ will increase so that this volume is not maximal, a contradiction. The inclusion
$Q_n\subset -nS$ means that  $-Q_n$ is a subset of a translate of
$nS$. Since  $-Q_n=Q_n$, we have $\alpha(S)\leq n$. But, thanks to \eqref{xi_S_alpha_S_geq_n}, $\alpha(Q_n;S)\geq n$ proving that $\alpha(Q_n;S)=n$. Let us also note that $d_i(S)\le 1$ because $S\subset Q_n$. From these properties and \eqref{alpha_d_i_S_eq} we have $d_i(S)=1$ for all  $i=1,...,n,$ proving the statement.
\hfill$\Box$

\smallskip
Note that this property of a maximum volume simplex was essentially used in \cite{nev_ukh_matzam SVFU_2019}.

 \smallskip
{\bf Corollary 2.7.} {\it For any $n$, we have $\xi_n(Q_n)\geq n.$}

\smallskip
This follows immediately from Corollary 2.3.
The simplicity of obtaining this estimate is apparent, since both the proofs of formula \eqref{alpha_d_i_S_eq} given in \cite{nevskii_dcg_2011} are quite labor-intensive.

 \smallskip
{\bf Corollary 2.8.} {\it For any $n$,
\begin{equation}\label{theta_n_3_minus_frac_ineq}
\theta_n(Q_n)\geq 3-\frac{4}{n+1}.
\end{equation}
Furthermore, if there exists a projector $P:C(Q_n)\to\Pi_1\left({\mathbb R}^n\right)$ satisfying the equality
\begin{equation}\label{P_on_cube_n_3_minus_frac_eq}
\|P\|_{Q_n}= 3-\frac{4}{n+1}
\end{equation}
then $P$ is minimal. Moreover, in this case $\xi_n(Q_n)=n.$}

\smallskip
{\it Proof.} Let us combine the right-hand inequality in \eqref{nev_ksi_P_ineq}  with the estimate
$\xi_n(Q_n)\geq n$:
\begin{equation}\label{nev_ksi_n_Q_n_theta_n_K_ineq}
n\leq\xi_n(Q_n) \leq
\frac{n+1}{2}\Bigl( \theta_n(Q_n)-1\Bigr)+1.
\end{equation}
Therefore, for any $n$, we have \eqref{theta_n_3_minus_frac_ineq}.  The equality
\eqref{P_on_cube_n_3_minus_frac_eq} means that we have an equality
also in \eqref{theta_n_3_minus_frac_ineq}. This means that the outer parts of \eqref{nev_ksi_n_Q_n_theta_n_K_ineq} are equal to $n$. Hence, in this case also
$\xi_n(Q_n)= n$.
\hfill$\Box$

\msk
Let us present a series of results devoted to calculation of the constant $\xi_n(Q_n)$ for small $n$.

The case $n=1$ is trivial: $\xi_1(Q_1)=1$ and a unique extremal simplex is the segment $[0,1]$ coinciding with $Q_1$.
Starting from $n=2$, the problem of calculating $\xi_n(Q_n)$  is rather difficult, especially, if
$n+1$ is not an Hadamard number. The values $\xi_2(Q_2)$ and  $\xi_3(Q_3)$ were discovered in \cite{nevskii_mais_2006_13_2}; see also
\cite{nevskii_monograph}.

For $n=2,$
$$\xi_2(Q_2)=\frac{3\sqrt{5}}{5}+1= 2.3416...$$
Up to rotations, the only extremal simplex is the triangle with vertices $(0,0),$ $(1,\tau),$
$(\tau,1)$, where
$\tau=(3-\sqrt{5})/2=0.3819\ldots$ .
 This number
satisfies
$\tau^2-3\tau+1=0$ or
$$\frac{\tau}{1-\tau}=  \frac{1-\tau}{1}.$$
Therefore,  $\tau$ delivers the so called {\it golden section} of the segment [0,1].
Sharp inequality
$\xi(Q_2;S)\geq 3\sqrt{5}/5+1$
for simplices $S\subset Q_2$ gives
new characterization of this classical notion.
See Fig. 1.

\begin{figure}[h!]
 \center{\includegraphics[scale=0.5]{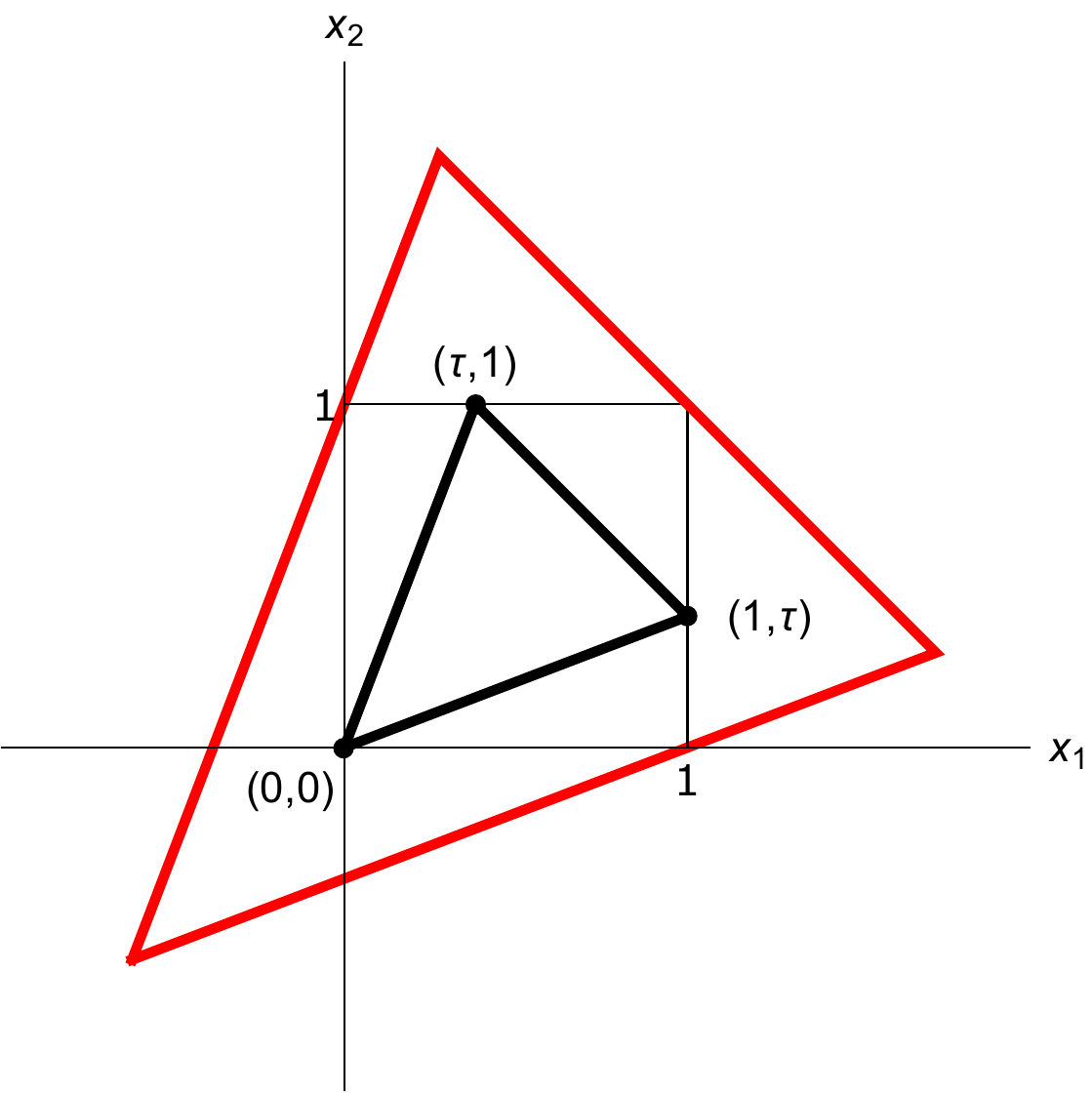}}
\caption{The case $n=2$.  The golden section simplex}
\label{fig:nev_ukl_simplex_n2}
\end{figure}


Combining our approach with some
results of \cite{lassak_dcg_1999}, we prove in \cite{nevskii_mais_2006_13_2} that  $\xi_3(Q_3)=3.$
Up to a~change of coordinates, each extremal simplex in $Q_3$ coincides with
the tetrahedron $S^\prime$ with vertices
$$(1,1,0),~~~(1,0,1),~~~(0,1,1),~~~(0,0,0)$$ or
the tetrahedron $S^{\prime\prime}$ with vertices
$$
\left(\frac{1}{2},0,0\right),~~~
\left(\frac{1}{2},1,0\right),~~~
\left(0,\frac{1}{2},1\right),~~~
\left(1,\frac{1}{2},1\right).
$$

In other words, if $\xi(Q_3;S)=3,$ then  either vertices of $S$ coincide with vertices of the cube and~form
a regular tetrahedron or coincide with the centers of opposite edges of two opposite faces of the cube and
does not belong to a common plane. One can see that all  axial diameters of $S^\prime$ and $S^{\prime\prime}$ are
equal to~1. This also follows from the equality $\xi_3(Q_3)=3.$ Note that $\vo(S^\prime)=\nu_3=
1/3,$
and $\vo(S^{\prime\prime})=
1/6$.

In \cite{nev_ukh_mais_2017_24_5} and  \cite{nev_ukh_posobie_2020} we show that each extremal simplex for $n=2$ and $n=3$ is {\it equisecting}, i.e., the~hyperplanes containing its faces cut off from the cube in outer direction from the simplex the~domains of equal volumes.

The estimate $\xi_n(Q_n)\geq n$ occurs to be exact in order of $n$. Indeed, let $n>2$ and let  $S^*$ be the simplex with the zero-vertex and other $n$ vertices coinciding with the  vertices of $Q_n$ adjacent to $(1,\ldots,1)$, i.e., with the vertices
\begin{equation}\label{vert_of_rigid_simplex}
(0,1,\ldots,1),~~~(1,0,\ldots,1),~~~\ldots,~~~(1,1,\ldots,0),~~~(0,0,\ldots,0)
\end{equation}
(see \cite{nevskii_mais_2011_2}, \cite[\S\,3.2]{nevskii_monograph}).
If  $n\geq 3$, then  $S^*$ has the following property
(see \cite[Lemma~3.3]{hudelson_1996}): replacement of an arbitrary vertex of  $S^*$ by
any point of $Q_n$ decreases the volume of the simplex.
For $n=2,3,4$ (and only in these cases),
simplex $S^*$ has maximum volume in $Q_n$ (see also \cite{hudelson_1996}).
 As it is shown in  \cite{nevskii_mais_2011_2}, if $n\geq 2$, then
$d_i(S^*)=1$ for every $i=1,...,n$, and therefore, thanks to \eqref{alpha_d_i_S_eq}, $\alpha(Q_n;S^*)=n.$
For this simplex,
$$\xi(Q_n;S^*)=\frac{n^2-3}{n-1}.$$

Hence, for $n>2$,
\begin{equation}\label{xi_n_leq_frac_general}
\xi_n(Q_n)\leq \frac{n^2-3}{n-1}.
\end{equation}
If $n>1$, the right-hand side of  \eqref{xi_n_leq_frac_general} is strictly smaller than  $n+1$.
Since $$
\xi_1(Q_1)=1,~~~\xi_2(Q_2)=\frac{3\sqrt{5}}{5}+1=2.3416\ldots,
$$
inequality $\xi_n<n+1$
holds true also for $n=1,2$.

We note that if the volume of simplex $S\subset Q_n$
is  maximal, then  $\xi(Q_n;S)\leq n+2$, see \cite{nevskii_monograph}.

Simplices satisfying the inclusions $S\subset Q_n\subset nS$ were studied in \cite{nev_ukh_mais_2017_24_5}.
The existence of such
simplices for a given $n$ is equivalent to the equality  $\xi_n(Q_n)=n.$  This holds not for each $n$,
but  nowaday a unique known dimension for which
$\xi_n(Q_n)>n$ is $n=2$.

Thus, always $n\leq \xi_n(Q_n)< n+1$, i.e., $\xi_n(Q_n)-n\in [0,1)$.
However,  the exact values of the constant
$\xi_n(Q_n)$ are currently only known for
$n=2$, $n=5$ and $n=9$, as well as for an infinite set of $n$ for which there exists an Hadamard matrix of order $n+1$.

In all these cases, except $n=2$, the equality $\xi_n(Q_n)=n$ holds.
In the noted Hadamard case, one can give the proof using the structure of Hadamard matrix of order $n+1$, see \cite{nev_ukh_beitrage_2018}.
In \cite{nev_ukh_beitrage_2018}, we have also discovered the exact values of $\xi_n(Q_n)$ for $n=5$ and $n=9$ and  constructed several infinite families of extremal simplices for $n=5,7,9$.

\begin{figure}[h!]
\center{\includegraphics[scale=0.4]{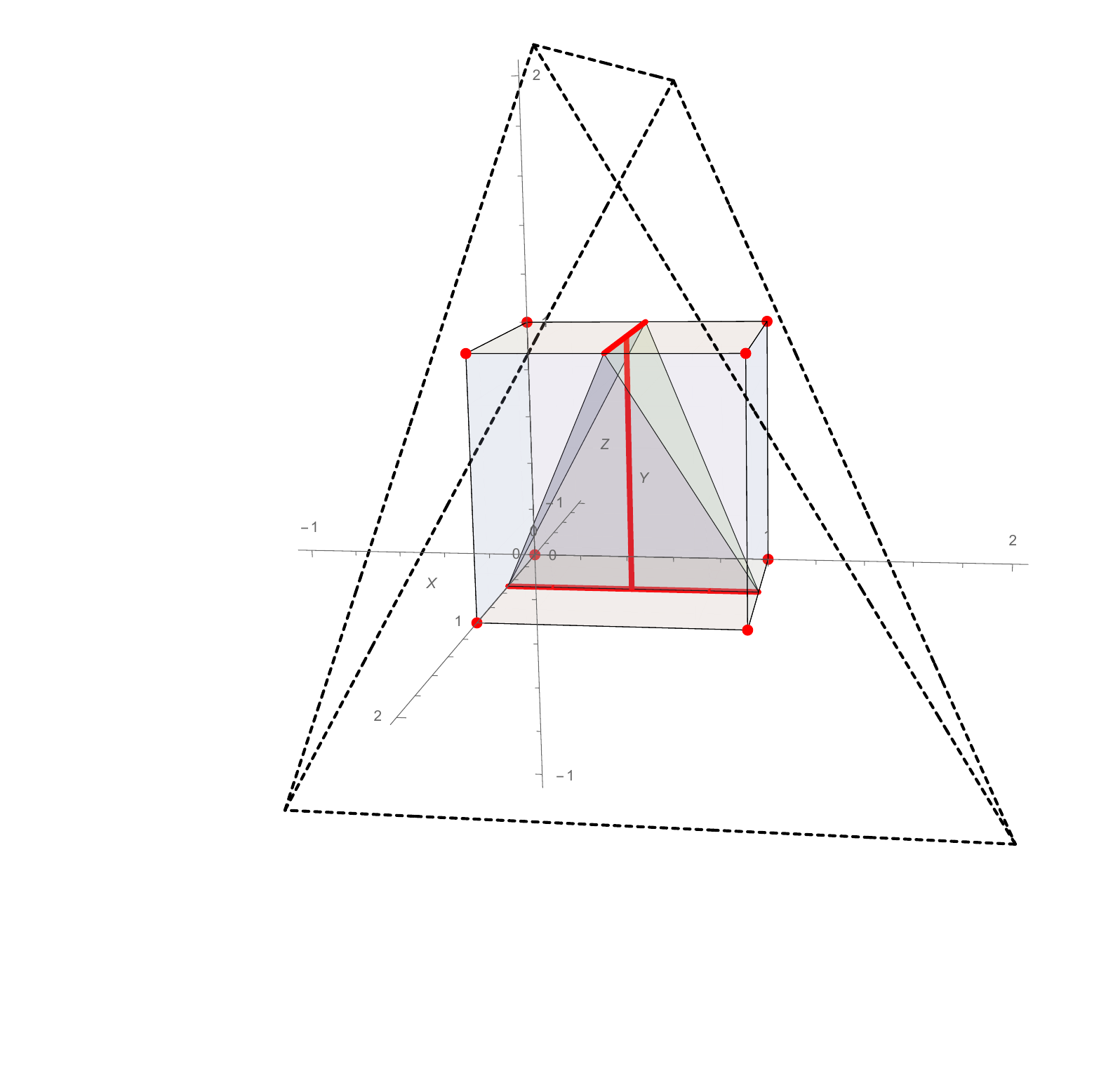}}
\vspace*{-20mm}
\caption{The case $n=3$. Simplex $S^{\prime\prime}$}
\label{fig:nev_ukl_simplex_s2}
\end{figure}

We call a nondegenerate simplex $S\subset Q_n$ {\it  perfect} provided $\xi(Q_n;S)=\xi_n(Q_n)$ and
the cube $Q_n$
is~inscribed into the simplex $\xi(Q_n;S)\,S$. For a long time, only two dimensions were known, namely $n=1$ and $n=3$, for which such simplices exist. Every three-dimensional perfect simplex is similar to $S^{\prime\prime}$.
In Fig. \ref{fig:nev_ukl_simplex_s2} it can be seen that
all vertices of $Q_3$ belong to the boundary of the simplex $3S^{\prime\prime}$.

\newpage

\begin{figure}[h!]
\center{\includegraphics[scale=0.3]{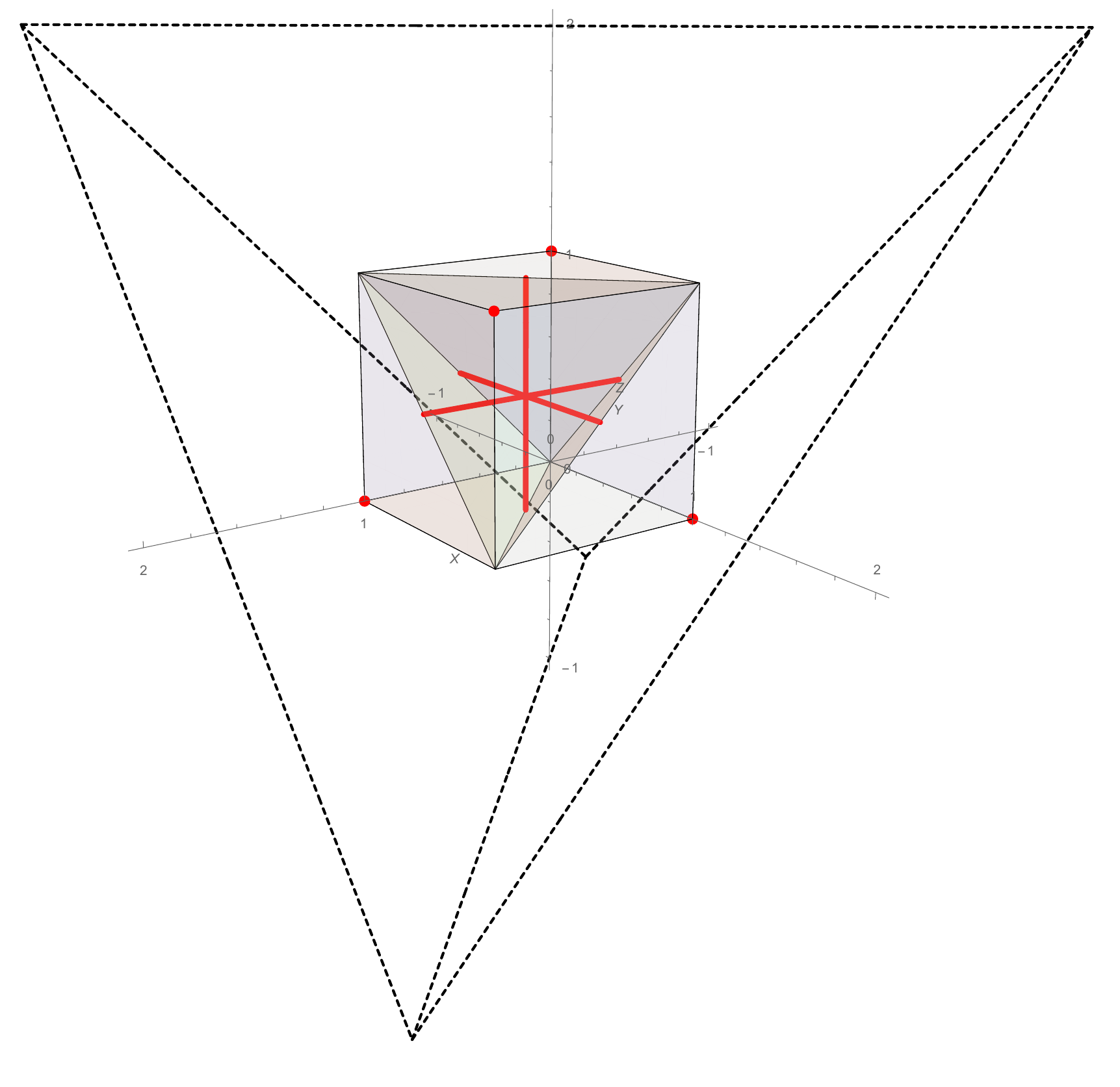}}
\caption{The case $n=3$. Simplex $S^\prime$}
\label{fig:nev_ukl_simplex_s1}
\end{figure}

On the contrary,
simplex  $S^{\prime}$ is not perfect, While $Q_3\subset 3S^\prime$, only four vertices of the cube
lie on~the~boundary of the
simplex $3S^{\prime}$. See Fig. \ref{fig:nev_ukl_simplex_s1}.

It was shown in \cite{nev_ukh_beitrage_2018} that perfect simplices do exist also in ${\mathbb R}^5$.

As we have mentioned above, nowaday $n=2$ is a unique known dimension for which $\xi_n(Q_n)>n.$

Let us also note that dimension $n=2$ remains the only even positive integer where we know the~sharp value of the constant $\xi_n(Q_n)$. Also, it is still unknown whether there exists an even $n$
such~that $\xi_n(Q_n)=n.$ The best of  known today upper bounds for $n=4$ and $n=6$ are
$$
\xi_4(Q_4)\leq \frac{19+5\sqrt{13}}{9}=4.1141\ldots, \quad \xi_6(Q_6)<6.0166.
$$
(see \cite{nev_ukh_mais_2016_23_5}, \cite{nev_ukh_mais_2017_24_1}).
 Using simplices of maximum volume in the cube, A.Yu. Ukhalov perfomed computer calculations to get as  accurate upper bounds for
$\xi_n(Q_n)$ as possible; the results for $n\leq 118$ are given in \cite{nev_ukh_mais_2018_25_1}
and also in the survey paper \cite{nevskii_ukh_sb_2021}.

Let us note that inequality $\xi_n(K)\leq n+2$ holds {\it for an arbitrary convex body $K\subset
{\mathbb R}^n$}. This is immediate from the following proposition proven by M. Lassak \cite{lassak_beitr_2011}: {\it If $S$ is a maximum volume simplex contained in $K$, then $\xi(K;S)\leq n+2.$} In other words, we have
\begin{equation}\label{S_sub_K_sub_n_plus_2_S}
S\subset K\subset (n+2)S.
\end{equation}

Note that this statement easily follows from formula \eqref{xi_K_S_formula}. In fact, because $S\subset K$ has
maximum volume, we have $|\Delta_j(x)|\leq |\Delta|$ for any $j=1,\ldots,n+1$ and  $x\in K.$
It suffices to use the relation \eqref{vol_S_det_A_eq} between determinants and volumes of simplices. 
Thanks to \eqref{Lagr_pol_thru_dets},
$$
-\lambda_j(x)\leq |\lambda(x)|=\frac{|\Delta_j(x)|}{|\Delta|} \leq 1, \quad x\in K.
$$
Thus, by \eqref{xi_K_S_formula}, we obtain
$$\xi(K;S)=(n+1)\max_{1\leq k\leq n+1}
\max_{x\in K}(-\lambda_k(x))+1\leq n+2.$$
 This approach was
used by the author in \cite{nevskii_monograph} in the case $K=Q_n$.

\SECT{3. The values of $\alpha(B_n;S)$ and $\xi(B_n;S)$}{3}
\label{alpha_ksi_ball}
\addtocontents{toc}{3. The values of $\alpha(B_n;S)$ and $\xi(B_n;S)$\hfill \thepage\par\VST}

\indent
\par In this section we discuss characteristics associated with the absorption of a Euclidean ball by~a~homothetic image of a simplex (with or without translation).
Replacing a cube with a ball makes many
questions much more simple. However, geometric interpretation of general results
has a certain interest also in this particular case.
Besides, we will note some new applications of the basic Lagrange polynomials.

Given an $n$-dimensional simplex $S$, let us introduce the following geometric characteristics. \linebreak {\it  The~ inradius of $S$} is the maximum of the radii of balls contained within $S$. The center of this unique maximum ball is called {\it the incenter of $S$.}
The boundary of the maximum ball is a sphere that has a~single common point
with each $(n-1)$-dimensional face of $S$. By {\it the circumradius of S}
we mean the minimum of the radii of balls containing $S$.
Note that the boundary of this unique minimal ball contains all the vertices of $S$ if and only if its center lies inside the simplex.

The inradius $r$ and the circumradius $R$ of a simplex $S$
satisfy the so-called {\it Euler inequality}
\begin{equation}\label{euler_ineq}
R\geq nr.
\end{equation}
Equality in
(\ref{euler_ineq}) takes place if and only if
$S$ is a regular simplex.
For the proof of the  Euler inequality, its history and various generalizations we refer the reader to \cite{klamkin_1979}, \cite{yang_wang_1985} and \cite{vince_2008}.

Let us note that the Euler inequality is equivalent to
the following statement:
{\it Suppose $B$ is a ball with radius $r_1$ and
$B^*$ is a ball with radius $r_2$. If
$B\subset S\subset B^*$, then $r_2\geq nr_1.$
Equality takes place if and only if $S$ is a regular simplex
inscribed into $B^*$ and
$B$ is the ball inscribed into $S$.} Another
equivalent form of these propositions is given by
Theorem 3.3.

Also we remark that the analog of the above property for parallelotopes is expressed by
Corollary~2.5.

We turn to presentation of some computational formulae for the quantities $\alpha(B_n;S)$ and $\xi(B_n;S)$ proven in \cite{nevskii_mais_2018_25_6}.
Let $x^{(1)},$ $\ldots,$
$x^{(n+1)}$ be the vertices and $\lambda_1,$ $\ldots,$
$\lambda_{n+1}$ be the basic Lagrange polynomials of a nondegenerate simplex
$S\subset {\mathbb R}^n$. As usual, $l_{ij}$ are the coefficients of $\lambda_j.$
See formulas from Section~1.

Given $j\in\{1,...,n+1\}$, let $\Gamma_j$ be the $(n-1)$-dimensional hyperplane determined by the equation $\lambda_j(x)=0$, and let
$\Sigma_j$ be the $(n-1)$-dimensional face of $S$ contained
in $\Gamma_j$. We let $h_j$ denote the~height of $S$ conducted from $x^{(j)}$ onto $\Gamma_j$, and we let $r$ denote the inradius of $S$. By $\sigma_j$ we denote $(n-1)$-dimensional measure of $\Sigma_j$. Finally, we put
$\sigma=\sum\limits_{j=1}^{n+1} \sigma_j$.

In \cite{nevskii_mais_2018_25_6} we show that the value $\alpha(B_n;S)$ can be calculated in various ways.

\smallskip
{\bf Theorem 3.1.} {\it The following equalities hold:
\begin{equation}\label{alpha_bs_sum_l_ij_equality}
\alpha(B_n;S)=
\sum_{j=1}^{n+1}\left(\sum_{i=1}^n l_{ij}^2\right)^{1/2},
\end{equation}
\begin{equation}\label{alpha_bs_h_j_equality}
\alpha(B_n;S)=\sum_{j=1}^{n+1}\frac{1}{h_j},
\end{equation}
\begin{equation}\label{alpha_bs_1_r_equality}
\alpha(B_n;S)= \frac{1}{r},
\end{equation}
$$\alpha(B_n;S)=\frac{\sigma}{n\,\vo(S)}.$$
}

\smallskip
From \eqref{alpha_bs_h_j_equality} and \eqref{alpha_bs_1_r_equality} we get
$$\frac{1}{r}=\sum_{j=1}^{n+1}\frac{1}{h_j}.$$
This geometric relation (which evidently can be obtained
in the direct way) occurs to be equivalent to~general formula \eqref{alpha_K_S_general_formula} for
$\alpha(K;S)$ in the
case when the convex body $K$ coincide with the  unit ball.

 It is interesting to compare
(\ref{alpha_bs_sum_l_ij_equality}) with the formula
$$\alpha(Q_n^\prime;S)=\sum_{i=1}^n\sum_{j=1}^{n+1} |l_{ij}|.$$
for  the cube $Q_n^\prime=[-1,1]^n$ (see Corollary 2.2).
Since $B_n\subset Q_n^\prime$, it follows that $\alpha(B_n;S)\leq \alpha(Q_n^\prime;S)$.
Of~course, analytically, this inequality is immediate from (\ref{alpha_bs_sum_l_ij_equality}) and the above formula for $\alpha(Q_n^\prime;S)$.

In turn, the inradius $r$ and the incenter $z$ of a simplex $S$ can be calculated as follows:
$$
r=\frac{1}{ \sum\limits_{j=1}^{n+1}\left(\sum\limits_{i=1}^n l_{ij}^2\right)^{1/2}}, \quad
z=\frac{1}{ \sum\limits_{j=1}^{n+1}\left(\sum\limits_{i=1}^n l_{ij}^2\right)^{1/2}}
\sum\limits_{j=1}^{n+1}\left(\sum\limits_{i=1}^n l_{ij}^2\right)^{1/2} x^{(j)}.
$$

For arbitrary $x^{(0)}$ and $\varrho>0$, the quantity
$\alpha\left(B(x^{(0)};\varrho);S\right)$ one can calculate
using the equality
$\alpha(B(x^{(0)};\varrho);S)$ $=$ $\varrho\alpha(B_n;S)$.

The general formula \eqref{xi_K_S_formula}
for  $\xi(K;S)$ in the case $K=B\left(x^{(0)};\varrho\right)$ provides the following (see \cite{nevskii_mais_2018_25_6}).

\smallskip
{\bf Theorem 3.2.} {\it
If $B\left(x^{(0)};\varrho\right)
 \not\subset S$, we have
\begin{equation}\label{ksi_b_x0_ro_s_l_ij_equality}
\xi\left(B\left(x^{(0)};\varrho\right);S\right)=
(n+1)\max_{1\leq j\leq n+1}
\left[\varrho\left(\sum_{i=1}^n l_{ij}^2\right)^{1/2}-
\sum_{i=1}^n l_{ij}x_i^{(0)}-l_{n+1,j}\right]+1.
\end{equation}
In particular, in the case
$B_n\not\subset S$
\begin{equation}\label{ksi_bs_l_ij_equality}
\xi(B_n;S)=
(n+1)\max_{1\leq j\leq n+1}\left[\left(\sum_{i=1}^n l_{ij}^2\right)^{1/2}
-l_{n+1,j}\right]+1.
\end{equation}
}

\smallskip
Now we turn to lower bounds for $\alpha(B_n;S)$ and $\xi(B_n;S)$ for $S\subset B_n.$
It is easy to show that the Euler inequality is equivalent to the following statement.

\smallskip
{\bf Theorem 3.3.}
{\it If $S\subset B_n$, then $\alpha(B_n;S)\geq n.$ The equality $\alpha(B_n;S)=n$ holds true if and only if
$S$ is~a~regular simplex inscribed into $B_n$.}

\smallskip
Recall  that the minimum value of
$\alpha(Q_n;S)$ for $S\subset Q_n$ is also equal to
$n$ (see \eqref{xi_S_alpha_S_geq_n}).
This value corresponds to those and only those
$S\subset Q_n$ for which every axial diameter
$d_i(S)$ is equal to $1$.
The~latter property is fulfilled for
the maximum volume simplices in $Q_n$, but, if $n>2$, not for the only ones.

If $S\subset B_n$, then $\xi(B_n;S)\geq n.$ The equality
$\xi(B_n;S)=n$ takes place if and only if
$S$ is a regular simplex inscribed into $B_n$.
These statements immediately follow from Theorem 3.3
and the inequa\-lity
$\xi(B_n;S)\geq \alpha(B_n;S)$. Paper \cite{nevskii_mais_2018_25_6} contains
also the direct arguments without applying the Euler inequality
that was used in \cite{nevskii_mais_2018_25_6} to obtain the estimate
$\alpha(B_n;S)\geq n$.

Now, let us consider an analog of the quantity
$\xi_n(Q_n)$, the quantity $\xi_n(B_n)$. In accordance with Definition 1.3,
$$
\xi_n(B_n)=\min \{ \xi(B_n;S): \,
S  \mbox{ is an $n$-dimensional simplex,} \,
S\subset B_n, \, \vo(S)\ne 0\}.
$$

While many problems concerning $\xi_n(Q_n)$ yet have not been solved (see Section 2),
the problem on~numbers $\xi_n(B_n) $ turns out to be simple. For any $n$, we have $\xi_n(B_n)=n$.
The only simplex $S\subset B_n$  extremal with respect to $\xi_n(B_n)$
is an arbitrary regular simplex inscribed into $B_n$.

Let us note that, thanks to \eqref{nev_ksi_P_ineq}, for every interpolation projector $P:C(B_n)\to\Pi_1\left({\mathbb R}^n\right)$ and the~corresponding
simplex $S\subset B_n$, the following inequalities
\begin{equation}\label{nev_ksi_P_ineq_ball}
\frac{n+1}{2n}\Bigl( \|P\|_{B_n}-1\Bigr)+1\leq
\xi(B_n;S) \leq
\frac{n+1}{2}\Bigl( \|P\|_{B_n}-1\Bigr)+1
\end{equation}
hold. Furthermore, if there exists an $1$-point in $B_n$ with respect to $S$
(see Section 1),
then the right-hand inequality in (\ref{nev_ksi_P_ineq_ball}) becomes an equality.
Hence,  
\begin{equation}\label{nev_ksi_n_B_theta_n_B_ineq}
\frac{n+1}{2n}\Bigl( \theta_n(B_n)-1\Bigr)+1\leq
\xi_n(B_n)=n \leq
\frac{n+1}{2}\Bigl(\theta_n(B_n)-1\Bigr)+1.
\end{equation}

Starting from $n=5$, the right-hand inequality in \eqref{nev_ksi_P_ineq_ball} is strict, while for $1\leq n\leq 4$ it turns into equality.
This fact is noted in Section 7.  Here is another argument, but also using some of the properties below.
Since we know the precise value
of $\theta_n(B_n)$ (see Section 9), the cases $1\leq n\leq 5$ are cheking directly. The result for $n\geq 6$
follows from the estimate $\theta_n(B_n)\geq \sqrt{n}$, since then
$$n<\frac{n+1}{2}\Bigl(\sqrt{n}-1\Bigr)+1.$$


\SECT{4. Estimates for $\theta_n(Q_n)$}{4}
\label{est_theta}
\addtocontents{toc}{4. Estimates for $\theta_n(Q_n)$\hfill \thepage\par\VST}

\indent\par Despite the apparent simplicity of formulation,
 the problem of finding exact values of $\theta_n(Q_n)$ is~very difficult.
Since 2006, these values are known only for $n=1,2,3,$ and $7$. Namely,
$$\theta_1(Q_1)=1, \quad \theta_2(Q_2)=\frac{2\sqrt{5}}{5}+1=1.8944\ldots, \quad \theta_3(Q_3)=2, \quad \theta_7(Q_7)= \frac{5}{2}.$$
Note that for each $n=1,2,3,7$ the right-hand relation in  \eqref{nev_ksi_n_K_theta_n_K_ineq}
becomes an equality:
\begin{equation}\label{righthand_eq_theta_ksi_cube}
\xi_n(Q_n) =
\frac{n+1}{2}\left( \theta_n(Q_n)-1\right)+1.
\end{equation}
If  $n=1, 3,$ or $7$, then an equality takes place also in \eqref{theta_n_3_minus_frac_ineq}, i.\,e.,
$$\theta_n(Q_n)= 3-\frac{4}{n+1}.$$
For $1\leq n\leq 3$, the  simplices corresponding to minimal projectors are just the same that the extremal simplices
with respect to $\xi_n(Q_n)$, see Section 2. In the seven-dimensional case, the nodes of minimal projector appear
from a unique (up to equivalence) Hadamard matrix of order $8$. The proofs are given in original paper
\cite{nevskii_mais_2009_16_1}   and also in \cite{nevskii_monograph} and \cite{nev_ukh_posobie_2022}.

Clearly, $\theta_1(Q_1)=1$. This value can be found also from \eqref{nev_ksi_n_K_theta_n_K_ineq}.
 For any projector,
there is a $1$-vertex $Q_1$ relative to the corresponding
simplex (which  is a
segment). Hence, we have \eqref{righthand_eq_theta_ksi_cube}. Now $\xi_1(Q_1)=1$ implies $\theta_1(Q_1)=1.$
 The nodes of a unique minimal  projector are 0 and 1.
 For the most complicated cases $n=2,3,$ see \cite{nevskii_monograph}.

Let us turn to the case $n=7$. Since $8$ is an Hadamard number, there exists a seven-dimensional regular simplex having vertices at vertices of the cube.
We can take the simplex with the vertices
$$
(1,1,1,1,1,1,1),~~~(0,1,0,1,0,1,0),~~~(0,0,1,1,0,0,1),~~~
(1,0,0,1,1,0,0),$$
$$
(0,0,0,0,1,1,1),~~~(1,0,1,0,0,1,0),~~~(1,1,0,0,0,0,1),~~~
(0,1,1,0,1,0,0).
$$

The equality $\xi_7(Q_7) =7$ (see Section 2) and
inequality
$$
\xi_n(Q_n) \leq \frac{n+1}{2}\left( \theta_n(Q_n)-1\right)+1
$$
imply $\theta_7(Q_7) \geq  \frac{5}{2}$. But for the corresponding projector, $\|P\|_{Q_7}=\frac{5}{2}$.
Therefore, $\theta_7 = \frac{5}{2}$, and this projector is minimal.

Calculations seems to be more easy if we make use of similarity reasons. Let us take the cube $Q_7^\prime=[-1,1]$
and consider the  projector
$P:C(Q_7^\prime)\to \Pi_1\left({\mathbb R}^7\right)$ constructed from  the Hadamard matrix
$$
{\bf H}_8=
\left( \begin{array}{cccccccc}
 1&1&1&1&1&1&1&1\\
 -1&1&-1&1&-1&1&-1&1\\
 -1&-1&1&1&-1&-1&1&1\\
 1&-1&-1&1&1&-1&-1&1\\
 -1&-1&-1&-1&1&1&1&1\\
 1&-1&1&-1&-1&1&-1&1\\
 1&1&-1&-1&-1&-1&1&1\\
 -1&1&1&-1&1&-1&-1&1
\end{array}\right).$$
(This matrix is unique up to the equivalence of Hadamard matrix of order 8.) Suppose the nodes of~$P$ are written in the rows of
this matrix except the last
column. Since
${\bf H}_8^{-1}$ $=$
$({1}/{8}){\bf H}_8^T,$  coefficients of $\lambda_j$ form the  $j$th row of
$({1}/{8}){\bf H}_8:$
$$
\lambda_1(x)=\frac{1}{8}
\left(x_1+x_2+x_3+x_4+x_5+x_6+x_7+1\right),
$$
$$\lambda_2(x)=\frac{1}{8}\left(-x_1+x_2-x_3+x_4-x_5+x_6-x_7+1\right),$$
and so on. Using \eqref{norm_P_cube_formula}, we get $\|P\|_{Q_7^\prime}={5}/{2}$, and, thanks to the above approach, $\theta_7(Q^\prime_7)={5}/{2}.$

However, $n=7$  is the biggest $n$ such that (i) $n+1$
is an Hadamard number and (ii) \eqref{P_on_cube_n_3_minus_frac_eq} holds provided $P$ corresponds to a regular simplex inscribed into the cube.
\msk

Now let us note some upper bounds for $\theta_n(Q_n)$.
For $n\geq 3$, the projector with  nodes \eqref{vert_of_rigid_simplex} satisfies
$\|P\|_{Q_n}\leq ({n+1})/{2}$ with the equality for odd $n$ (\cite{nevskii_mais_2003_10_1},
 \cite{nevskii_monograph}). Therefore, if $n\ne 2$, then 
 $$\theta_n(Q_n)\leq 
\frac{n+1}{2}.$$

Recall that
$h_n$ denotes the maximum value of a determinant of order $n$ with entries $0$ or~$1$
(see~Definition 1.9).  Suppose $P$ corresponds to a maximum volume simplex $S\subset Q_n$  with some vertex coinciding with a vertex of the cube. Then
\begin{equation}\label{th1_formula}
\|P\|_{Q_n}\leq \frac{2h_{n+1}}{h_n}+1.
\end{equation}
Hence,
$$
\theta_n(Q_n)\leq \frac{2h_{n+1}}{h_n}+1
$$
for all $n$ and also $\theta_n(Q_n)\leq \sqrt{2n+3}+1$,
provided $n+1$ is an Hadamard number (see \cite{nevskii_mais_22}).
Earlier (see \cite{nevskii_mais_2003_10_1}  and \cite{nevskii_monograph})  it was proved that in the Hadamard case  $\xi_n(Q_n)\leq
\sqrt{n+1}$. Let us give here the~proof from \cite{nevskii_mais_22} essentially making use of the structure of an Hadamard matrix.

\smallskip
{\bf Theorem 4.1.} {\it  Let $n+1$ be an Hadamard number, and let $S$ be an $n$-dimensional regular simplex having the vertices at vertices of  $Q_n^\prime$. Then, for the
corresponding interpolation projector $P:C(Q_n^\prime)\to \Pi_1({\mathbb R}^n)$, we have
\begin{equation}\label{th2_formula}
\|P\|_{Q_n^\prime}\leq \sqrt{n+1}.
\end{equation}
}

\smallskip
{\it Proof.}  Let  $x^{(1)}, \ldots,  x^{(n+1)}$ be the vertices and  let $\lambda_1, \ldots,  \lambda_{n+1}$ be the basic Lagrange polynomials of~the~simplex.
Thanks to the theorem's hypothesis, the vertex matrix $\bf S$ is an~Hadamard matrix of~order  $n+1$ with the last column consisting of $1$'s.

Let us show that
\begin{equation}\label{proof_th2_1}
\sum_{j=1}^{n+1} \lambda_j(x)^2=\frac{||x||^2+1}{n+1}
~~~\text{for every}~~~x\in{\mathbb R}^n.
\end{equation}

Let $y=(x_1,\ldots,x_n,1)\in {\mathbb R}^{n+1}.$
The coefficients of
$\lambda_j$ form the $j$th column of~${\bf S}^{-1}$. Since  $\bf S$ is an Hadamard matrix, it satisfies~\eqref{bfS_is_Hadamard}, whence
$$\lambda_j(x)=\frac{1}{n+1}\langle h^{(j)},y\rangle.$$
Further,  
${h^{(j)}}/{\sqrt{n+1} }$  form an ortho\-normalized basis of ${\mathbb R}^{n+1},$ therefore,
$$y=\sum_{j=1}^m \frac{ \langle h^{(j)},y\rangle }{\sqrt{n+1}}
\frac{h^{(j)}}{\sqrt{n+1}}, \quad \langle y, y\rangle=\sum_{j=1}^{n+1}
\frac{\langle h^{(j)},y\rangle^2} {n+1}.$$
From this, we have
$$||x||^2+1=\langle y, y\rangle =\sum_{j=1}^{n+1}
\frac{\langle h^{(j)},y\rangle^2} {n+1}=(n+1)\sum_{j=1}^{n+1} \lambda_j(x)^2$$
proving \eqref{proof_th2_1}.

If $x$ is a vertex of $Q_n^\prime$, then  $||x||^2=n$ and \eqref{proof_th2_1} implies
$\sum \lambda_j(x)^2=1.$
Applying Cauchy inequality, for each $x\in \ver(Q_n^\prime)$, we have
$$\sum_{j=1}^{n+1}|\lambda_j(x)|\leq \left(\sum_{j=1}^{n+1} \lambda_j(x)^2\right)^{\frac{1}{2}}
\cdot\sqrt{n+1}=\sqrt{n+1}.$$

Now, let $P:C(Q_n^\prime)\to \Pi_1({\mathbb R}^n)$ be an interpolation projector corresponding to $S$.
Then, thanks to~(\ref{norm_P_cube_formula}),
$$
\|P\|_{Q_n^\prime}=\max_{x\in\ver(Q_n^\prime)} \sum_{j=1}^{n+1} |\lambda_j(x)|\leq \sqrt{n+1}
$$
completing the proof of the theorem.
\hfill$\Box$

\smallskip
By similarity reasons, Theorem 4.1 overcomes to any $n$-dimensional cube, e.\,g., to the cube
$Q_n=[0,1]^n$. Thus, if $n+1$ is an Hadamard number, then $\theta_n(Q_n)\leq \sqrt{n+1}.$

We can slightly improve this estimate by using interpolation on a ball.

Let $p_n=\max\{\psi(a),\psi(a+1)\}$ where
$$
\psi(t):=\frac{2\sqrt{n}}{n+1}\Bigl(t(n+1-t)\Bigr)^{1/2}+
\left|1-\frac{2t}{n+1}\right|, \quad
a=\left\lfloor\frac{n+1}{2}-\frac{\sqrt{n+1}}{2}\right\rfloor.
$$

 Now, let us assume that $n+1$ is an Hadamard number.
Consider an interpolation projector $P$ with~the~nodes at~those vertices of the cube $Q_n^\prime$ that form a regular simplex $S$. Since  $Q_n^\prime$ is inscribed into the unit
ball $B_n$,   simplex $S$ is also inscribed into $B_n$. It remains to apply  formula
\eqref{norm_P_intro_cepochka} for the~projector's norm both on the cube and on the ball:
$$\|P\|_{Q_n^\prime}=\max_{x\in Q_n^\prime}\sum_{j=1}^{n+1}
|\lambda_j(x)|\leq \max_{x\in B_n}\sum_{j=1}^{n+1}
|\lambda_j(x)|= \|P\|_{B_n}=p_n\leq \sqrt{n+1}.$$

(For the equality $\|P\|_{B_n}=p_n$ and the final estimate, we refer the reader to Section 7.) Hence,
$$
\theta_n(Q_n)\leq p_n\leq \sqrt{n+1}.
$$

The equality
 $\|P\|_{Q_n^\prime}=\sqrt{n+1}$ may hold as for all regular simplices having vertices at vertices of~the~cube $(n=1, n=3)$, as for some of them $(n=15),$
 or  may not be executed at all.

As it is shown in Section 6, for each $n$,
$\theta_n(Q_n)\geq \sqrt{n-1}/e.$ Therefore,
if $n+1$ is an Hadamard number, then
$$\frac{\sqrt{n-1}}{e}\leq \theta_n(Q_n)\leq  \sqrt{n+1}.$$

\msk
The upper bounds of $\theta_n(Q_n)$ for special $n$ were improved by A. Ukhalov and his students 
\linebreak with the~help of computer methods. In particular,  simplices of maximum volume in the cube were considered.
In all situations where $n+1$ is an Hadamard number, the full set of Hadamard matrices of the corresponding order was used. In particular, to obtain an estimate for $\theta_{23}$, all existing $60$ Hadamard matrices of order $24$ were considered. To estimate $\theta_{27}$,
we have to consider $487$ Hadamard matrices of order $28$.
Known nowaday upper estimates for $1\leq n\leq 27$ are given in \cite{nev_ukh_posobie_2022}.
(Here, for~brevity, $\theta_n=\theta_n(Q_n)$.)
$$\theta_1=1, \quad \theta_2=\frac{2\sqrt{5}}{5}+1, \quad  \theta_3=2,
\quad \theta_4\leq \frac{3  (4+\sqrt{2})}{7},\quad\theta_5\leq  2.448804,$$
$$\theta_6\leq 2.6000\ldots ,\quad \theta_7=\frac{5}{2}, \quad \theta_8\leq \frac{22}{7},\quad\theta_9\leq 3, \quad\theta_{10}\leq \frac{19}{5}, \quad\theta_{11}\leq 3,$$
$$\theta_{12}\leq  \frac{17}{5},\quad
\theta_{13}\leq  \frac{49}{13},\quad \theta_{14}\leq \frac{21}{5},
\quad \theta_{15}\leq \frac{7}{2}, \quad \theta_{16}\leq  \frac{21}{5},\quad \theta_{17}\leq \frac{139}{34},
$$
 $$ \theta_{18}\leq 5.1400\ldots,
\quad \theta_{19}\leq 4, \quad \theta_{20}\leq 4.68879\ldots, \quad \theta_{21}\leq
\frac{251}{50},
\quad\theta_{22}\leq  \frac{1817}{335},$$
$$\theta_{23}\leq  \frac{9}{2},
\quad \theta_{24}\leq  \frac{103}{21}, \quad \theta_{25}\leq  5,
\quad \theta_{26}\leq  \frac{474}{91},
\quad \theta_{27}\leq  5.$$

Note that $n$-dimensional regular simplices with vertices at the vertices of the cube can be located differently with respect to  the vertices and faces of the cube.
This is quite evident when  the norms of the corresponding projectors are different. But this is also possible if regular simplices generate the same norms.
In \cite{nevskii_mais_22} we describe an approach based on comparison of $\mu$-vertices of the cube with~respect to~various~simplices. 

Recall that a $\mu$-vertex of the cube $Q$ with respect to a simplex $S\subset Q$  is such
a vertex $x$ of~the~cube that for the corresponding projector $P=P_S$ we have
$\|P\|_Q=\sum |\lambda_j(x)|$ and exactly $\mu$ numbers $\lambda_j(x)$ are negative
(see Definition 1.7).
Of course, simplices which have different sets of $\mu$-vertices with respect to~the~containing cube, are differently  located in the cube, even if they have the same projector's norms.
These simplices are non-equivalent in the following sense: one of them cannot be mapped into another by an orthogonal transform which maps the cube into~itself.
Let us give some examples for the case $Q=Q_n^\prime.$

While obtaining estimates for minimal norms of projectors, in \cite{nev_ukh_mais_2018_25_3} various $n$-dimensional regular simplices arising from different Hadamard matrices of~order $n+1$ were calculated.
 For~$n=15$, the~order of matrices is equal to  $16$. Up~to~equivalence, there are exactly five  Hadamard matrices. They correspond to five simplices described in Table \ref{tab:nev_Reg_simplices_for_n_15}.

\begin{table}[!htbp]
\begin{center}
\caption{\label{tab:nev_Reg_simplices_for_n_15} Regular simplices for $n=15$}
\medskip
\bgroup
$
\def\arraystretch{1.5}
\begin{array}{|c|c|c|c|}
  \hline
 S & {\|P\|}_{Q_{15}^\prime} & \mbox{Values } \ \mu& \mbox{ Number of} \ \mu\mbox{-vertices}  \\
  \hline
 S_1 & 4 & 6& m_6=448 \\
  S_2 & 4& 6& m_6=192 \\ S_3 & 4& 6& m_6=64 \\ S_4 & \frac{7}{2} & 4, 5, 6, 8& m_4=896, m_5=1344, m_6=5376, m_8=1344 \\
   S_5 & \frac{7}{2} & 4, 5, 6, 8& m_4=896, m_5=1344, m_6=5376, m_8=1344 \\
  \hline
 \end{array}
$
\egroup
\end{center}
\end{table}

By $m_\mu$ we denote a number of $\mu$-vertices of the cube  $Q_n^\prime$ with respect to every simplex. For the rest
$1\leq\mu\leq 15$,
excepting the given in Table \ref{tab:nev_Reg_simplices_for_n_15},
values  $m_\mu$ are~zero.
Every simplex $S_1$, $S_2$, and~$S_3$ generates the same projector's norm and has only $6$-vertices. But the numbers of $6$-vertices for~them are different;
hence, these simplices are pairwise non-equivalent. Each of them also is non-equivalent  as to $S_4$, so as to $S_5.$
The latter simplices generates equal norms and have the same sets \linebreak of $\mu$-vertices. Also we have $\theta_{15}\leq {7}/{2}=3.5$.
This is more exact than the estimate  $\theta_{15}(Q_{15})\leq 4$  provided by \eqref{th2_formula}.

Another example given in \cite{nevskii_mais_22}  is related to $n=23.$
Regular inscribed simplices were built from the available 60~pairwise non-equivalent Hadamard matrices of order  $24$. For all simplices, excepting ones with the numbers 16, 53, 59, and 60,
the projector norm  equals ${14}/{3}=4.6666\ldots$, while for~each of these four simplices the projector norm is ${9}/{2}=4.5$. In~particular, we have  $\theta_{23}(Q_{23})\leq 4.5$
(which was also noted in  \cite{nev_ukh_mais_2018_25_3}). This inequality is more exact than the estimate  $\theta_{23}(Q_{23})\leq \sqrt{24}=4.8989\ldots$.  Each  of the four exceptional simplices is not equivalent to any of the 56 others.

Despite the possible differs, for all regular simplices with vertices at the vertices of the cube, inequality \eqref{th2_formula} holds. The inscribed regular simplices satisfying $\|P\|_{Q_n^\prime}=\sqrt{n+1}$ exist at least for~$n=1$, $n=3$, and $n=15.$ The problem of full description of dimensions $n$ with such property is still open.

The best  nowaday known  lower bound of $\theta_n(Q_n)$ has the form
 \begin{equation}\label{nev_theta_n_max_ineq}
\theta_n(Q_n)
\geq  \max \left[
3 - \frac{4}{n+1}, \,
\chi_n^{-1} \left(\frac{1}{\nu_n}\right) \right].
\end{equation}
Here $\chi_n$ is the standardized Legendre polynomial of degree $n$, see Section 5.

The values of the right-hand side of \eqref{nev_theta_n_max_ineq} for $1\leq n\leq 54$ are given in \cite{nev_ukh_posobie_2022}.
Inequality $$\chi_n^{-1}\left(\frac{1}{\nu_n}\right)> 3-\frac{4}{n+1}$$ takes place for $n\geq 53$.
On estimates via the function $\chi_n^{-1}$  see Section 6.

As noted in Section 1 (see the right-hand relation in \eqref{nev_ksi_P_ineq}),
the inequality
\begin{equation}\label{nev_ksi_n_teta_n_ineq_sec4}
\xi_n(Q_n) \leq
\frac{n+1}{2}\left( \theta_n(Q_n)-1\right)+1
\end{equation}
 is true for any $n$.
So far, we  know only four values of $n$ for which this relation becomes an
equality: $n=1, 2, 3$, and $7$. These are exactly the cases in which
we know the exact values of $\theta_n(Q_n)$ and $\xi_n(Q_n)$. In
\cite{nev_ukh_mais_2016_23_5}
the authors conjectured that the minimum of $n$ for which inequality (\ref{nev_ksi_n_teta_n_ineq_sec4}) is strict is $4.$ This is still an open problem.

We note that, thanks to the equivalence $\xi_n(Q_n)\asymp n$ and inequality $\theta_n\geq c \sqrt{n}$, for all sufficiently large $n$, we have
\begin{equation}\label{nev_strict}
\xi_n(Q_n) <
\frac{n+1}{2}\left(\theta_n(Q_n)-1\right)+1.
\end{equation}

Let $n_0$ be the minimal natural number
such that for all $n\geq n_0$ inequality (\ref{nev_strict}) holds. The problem about the exact value of $n_0$ is very difficult. The known
lower and upper bounds differ quite significantly.
From the preceding, we have the estimate $n_0\geq 8.$
In 2009 we proved that $n_0\leq 57$ (see \cite{nevskii_mais_2009_16_1},
\cite[\S\,3.7]{nevskii_monograph}).
A sufficient condition for the validity of
(\ref{nev_strict})
for $n>2$ is the inequality
\begin{equation}\label{nev_chi_nu_ineq}
\chi_n \left(\frac{3n-5}{n-1}\right)\cdot\nu_n<1.
\end{equation}

It was proved  in \cite{nevskii_mais_2009_16_1} that \eqref{nev_chi_nu_ineq}  is satisfied for $n\leq 57$.
Later calculations allowed the upper bound of $n_0$ to be slightly lowered. These results are described in \cite{nev_ukh_mais_2018_25_3} where  it is noted
that $n_0\leq 53$. In~other words,
inequality \eqref{nev_strict} is satisfied at least
starting from $n=53$. Note that a better estimate from above for $n_0$ is an open problem.


\SECT{5. Legendre polynomials and the measure of $E_{n,\gamma}$}{5}\label{legendre_pol_mes_E}
\addtocontents{toc}{5. Legendre polynomials and the measure of $E_{n,\gamma}$\hfill \thepage\par\VST}

\indent\par {\it The standardized Legendre polynomial of degree $n$} is the function
$$
\chi_n(t)=\frac{1}{2^nn!}\left[ (t^2-1)^n \right] ^{(n)}
$$
(Rodrigues' formula). For properties of $\chi_n$ see, e.g.,
\cite{suetin_1979},\cite{sege_1962}.
Legendre polynomials are orthogonal on the segment
$[-1,1]$
with the weight $w(t)\equiv1.$
The first  standardized  Legendre polynomials are
$$\chi_0(t)=1, \quad
\chi_1(t)=t, \quad
\chi_2(t)=\frac{1}{2}\left(3t^2-1\right), \quad
\chi_3(t)=\frac{1}{2}\left(5t^3-3t\right),$$
$$\chi_4(t)=\frac{1}{8}\left(35t^4-30t^2+3\right), \quad
\chi_5(t)=\frac{1}{8}\left(63t^5-70t^3+15t\right).$$
We have
\begin{equation}\label{reccurent_Legendre}
\chi_{n+1}(t)=
\frac{2n+1}{n+1}\,t\chi_n(t)
-\frac{n}{n+1}\,\chi_{n-1}(t).
\end{equation}
This implies
$\chi_n(1)=1$. If $n\geq 1$, then $\chi_n(t)$
increases as~$t\geq 1$. The latter properties also easily follow from \eqref{th_2_1_1}.

We let $\chi_n^{-1}$ denote the function inverse to $\chi_n$ on the semi-axis $[1,+\infty)$.

The appearance of Legendre polynomials in the circle of our questions 
is due to some their property.  For $\gamma\geq 1$, let us  define the set
$$E_{n,\gamma}:= \left\{ x\in{\mathbb R}^n :
\sum_{j=1}^n |x_j| +
\left|1- \sum_{j=1}^n x_j\right| \leq \gamma \right\}.$$

In \cite{nevskii_mais_2003_10_1} we proved such a statement (the proof is also given in \cite{nevskii_monograph} and \cite{nev_ukh_posobie_2022}).

\smallskip
{\bf Theorem 5.1.} {\it The  following equalities hold:
\begin{equation}\label{th_2_1_1}
\mes(E_{n,\gamma})
= \frac{1}{2^n n!} \sum_{i=0}^{n}  {n \choose i}^2 (\gamma - 1)^{n-i}
(\gamma + 1)^i = \frac{\chi_n(\gamma)}{n!}.
\end{equation}
}

\smallskip
{\it Proof.} First let us  establish the left-hand equality in \eqref{th_2_1_1}.
Let
$$
E^{(1)}=\{x\in E_{n,\gamma}: \sum x_i>1 \},~~~~ E^{(2)}=\{x\in E_{n,\gamma}:\sum x_i \leq 1\}.
$$

Let us find sequentially $m_1=\mes(E^{(1)})$ and $m_2=\mes(E^{(2)}).$

Temporarily fix $k$, $1\le k\le n$, and consider a non-empty subset
$G\subset E^{(1)}$ consisting of all points $x=(x_1,...,x_n)$ such that
$x_1,\ldots,x_k\geq 0$ and $x_{k+1},\ldots, x_n< 0$. Let $y_i=x_i$ for $i=1,...,k$,
$y_{i}=-x_{i}$ for~$i=k+1,...,n$, and let $y=(y_1,...,y_n)$. Then
$$G=\left\{y: 1+y_{k+1}+\ldots+y_n\leq y_1+\ldots+y_k\leq \frac{\gamma+1}{2}, \
y_i\geq 0 \right \},$$
hence
$$\mes(G)=
\int\limits_{1}^{\alpha} \,dy_1
\int\limits_{1}^{\alpha-y_1} \,dy_2
\ldots
\int\limits_{1}^{\alpha-y_1-\ldots-y_{k-1}} \,dy_k\cdot$$
$$\cdot
\int\limits_{0}^{y_1+\ldots+y_k-1} \,dy_{k+1}
\int\limits_{0}^{y_1+\ldots+y_k-1-y_{k+1}} \,dy_{k+2}
\ldots
\int\limits_{0}^{y_1+\ldots+y_k-1-y_{k+1}-\ldots-y_{n-1}} \,dy_n.$$
Throughout the proof, $\alpha=
(\gamma+1)/2.$
If $b>0$, then
$$\int\limits_{0}^{b} \,dz_1
\int\limits_{0}^{b-z_1} \,dz_2\ldots\int\limits_{0}^{b-z_1-\ldots-z_{l-1} }
\,dz_l=\frac{b^l}{l!},$$
so
$$\mes(G)=$$
$$=
\int\limits_{1}^{\alpha} \,dy_1
\int\limits_{1}^{\alpha-y_1} \,dy_2
\ldots
\int\limits_{1}^{\alpha-y_1-\ldots-y_{k-1}}
\frac{1}{(n-k)!}{\left(y_1+\ldots+y_k-1\right)}^{n-k}
\,dy_k=$$
$$=\left (\int\limits_{\ y_1+\ldots+y_k\leq\alpha} -\int\limits_{y_1+\ldots+y_k\leq 1}
\right)
\frac{1}{(n-k)!}  {(y_1+\ldots+y_k-1)}^{n-k}
\,dy_1\ldots \,dy_k=$$
$$= J_1 - J_2.$$
The first integral equals
$$J_1=
\sum_{j=1}^{k} (-1)^{j+1} \frac{(\alpha-1)^{n-k+j}}{(n-k+j)!}
\frac{\alpha^{k-j}}{(k-j)!}+\frac{(-1)^{n+k}}{n!}
.$$
The value of $J_2$ appears from this expression
if instead of $\alpha$ we take 1. Consequently,

$$\mes(G)= J_1 - J_2=
\sum_{j=1}^{k} (-1)^{j+1} \frac{ (\alpha-1)^{n-k+j} }{(n-k+j)!}
\frac{ \alpha^{k-j} }{ (k-j)! }
=$$
$$
=\frac{ (-1)^{k+1} }{n!}
\sum_{i=0}^{k-1} {n\choose i}
(\alpha-1)^{n-i}(-\alpha)^i.$$

Clearly, the set $E^{(1)}$ is the union of all pairwise disjoint sets $G$ with various
$k=1,\ldots,n,$. Therefore, the measure of $E^{(1)}$ is equal to
$$
m_1=\sum_{k=1}^{n} {n\choose k} \frac{(-1)^{k+1}}{n!} \sum_{i=0}^{k- 1}
{n \choose i} (\alpha-1)^{n-i}(-\alpha)^i.
$$

Changing the order of summation and using the identity
\begin{equation}\label{th_2_1_2}
\sum_{k=0}^i (-1)^k {n\choose k}
= (-1)^i {{n-1}\choose i}
\end{equation}
(see, e.\,g., \cite{prudnikov_1981}) we get
 \begin{equation}\label{th_2_1_3}
m_1=
\frac{1}{n!} \sum_{i=0}^{n-1} {n \choose i} (\alpha-1)^{n-i}(-\alpha)^i
\sum_{k=0}^i (-1)^k {n\choose k}
=\frac{1}{n!} \sum_{i=0}^{n-1} {n \choose i} { {n-1}\choose i}
  (\alpha-1)^{n-i}\alpha^i.
  \end{equation}

Now let us turn to $E^{(2)}.$ First, note that $E^{(2)}$
contains the domain
$S:=\{x_i\geq 0,$ $ \sum x_i \leq 1\},$ with the measure  
$1/n!$. Next,
fix $k\in\{1,...,n\}$ and consider the subset
$G^\prime \subset E^{(2)}$
corresponding to the inequalities $x_1,\ldots,x_k<0;$ $x_{k+1},\ldots,x_n\geq
0.$ Put $y_1=-x_1,\ldots,y_k=-x_k;$ $y_{k+1}=x_{k+1},\ldots,y_n=x_n.$ Then
$$
G^\prime=\{y: y_{k+1}+\ldots+y_n\leq 1+ y_1+\ldots+y_k\leq \frac{\gamma-1}{2}, \
y_i\geq 0 \}.
$$

Let $\beta=({\gamma-1})/{2}.$ Then the following equalities hold:
$$\mes(G^\prime)=
\int\limits_{0}^{\beta} \,dy_1
\int\limits_{0}^{\beta-y_1} \,dy_2
\ldots
\int\limits_{0}^{\beta-y_1-\ldots-y_{k-1}} \,dy_k\cdot$$
$$\cdot
\int\limits_{0}^{1+y_1+\ldots+y_k} \,dy_{k+1}
\int\limits_{0}^{1+y_1+\ldots+y_k-y_{k+1}} \,dy_{k+2}
\ldots
\int\limits_{0}^{1+y_1+\ldots+y_k-y_{k+1}-\ldots-y_{n-1}} \,dy_n=$$
$$=
\int\limits_{0}^{\beta} \,dy_1
\int\limits_{0}^{\beta-y_1} \,dy_2
\ldots
\int\limits_{0}^{\beta-y_1-\ldots-y_{k-1}} \frac{(1+y_1+\ldots+y_k)^{n-k}}{(n-k)!}
\,dy_k=$$
$$=\sum_{j=0}^{k-1} (-1)^{k-1-j} \frac{(1+\beta)^{n-j}\beta^j}{(n-j) !j!}
+\frac{(-1)^k}{n!}=$$
$$=\frac{(-1)^{k+1}}{n!}
\Bigl(\sum_{j=0}^{k-1} {n\choose j} (1+\beta)^{n-j}(-\beta)^j - 1 \Bigr).$$

Clearly, the set $E^{(2)}$ is the union of all such sets $G^\prime$
corresponding to various $k=1,\ldots,
n,$ and~the~simplex $S.$ Therefore,
$$
m_2=\mes(E^{(2)})=$$
$$=
\frac{1}{n!} \Bigl( \sum_{k=1}^{n} (-1)^{k+1} {n\choose k}
\Bigl(\sum_{j=0}^{k-1} {n\choose j} (1+\beta)^{n-j}(-\beta)^j - 1 \Bigr)
+1 \Bigr).$$
Remark that
$$\displaystyle 1+\beta=\frac{\gamma+1}{2}=\alpha, \quad
\beta=\frac{\gamma-1}{2}=\alpha-1.$$
Using the substitution $i=n-j$ in the internal sum we obtain
$$m_2= $$
$$=
\frac{1}{n!} \Bigl(1+\sum_{k=1}^{n} (-1)^{k+1} {n\choose k}
\Bigl( (-1)^n\sum_{i=n-k+1}^{n} {n\choose i} (\alpha-1)^{n-i}(-\alpha)^i
  - 1\Bigr)
\Bigr)=$$
$$=
\frac{(-1)^n}{n!} \sum_{k=1}^{n} (-1)^{k+1} {n\choose k}
\sum_{i=n-k+1}^{n} {n\choose i} (\alpha-1)^{n-i}(-\alpha)^i.$$
We took into account that
$$\sum_{k=0}^{n} (-1)^k {n\choose k} = \sum_{k=1}^{n} (-1)^k
{n\choose k} + 1 = 0.$$
Changing the order of summation, we obtain
$$m_2=
\frac{(-1)^n}{n!}
\sum_{i=1}^{n} {n\choose i} (\alpha-1)^{n-i}(-\alpha)^i
\sum_{k=n+1-i}^{n} (-1)^{k+1} {n\choose k}.$$
Using \eqref{th_2_1_2}, we can write
$$
\sum_{k=n+1-i}^{n} (-1)^{k+1} {n\choose k}
=\sum_{k=n+1-i}^{n} (-1)^{k+1}{n\choose {n-k}}
=\sum_{j=0}^{i-1} (-1)^{n-j+1}{n\choose j}=
(-1)^{n+i} { {n-1}\choose {i-1} }.$$
Thus we have
\begin{equation}\label{th_2_1_4}
m_2=
\frac{1}{n!}\sum_{i=1}^{n} {n\choose i} { {n-1}\choose {i-1} }
  (\alpha-1)^{n-i}
\alpha^i.
\end{equation}

Thanks to equalities \eqref{th_2_1_3} and \eqref{th_2_1_4}, we have
$$
\mes(E_{n,\gamma})=m_1+m_2=
$$
$$
=\frac{1}{n!} \sum_{i=0}^{n-1} {n \choose i} { {n-1}\choose i}
  (\alpha-1)^{n-i}\alpha^i +
\frac{1}{n!}\sum_{i=1}^{n} { n\choose i} { {n-1}\choose {i-1} }
  (\alpha-1)^{n-i}
\alpha^i=$$
$$=\frac{1}{n!} \sum_{i=1}^{n-1} { n \choose i}
\left( { {n-1}\choose i} + { {n-1}\choose {i-1} } \right)
  (\alpha-1)^{n-i}\alpha^i +$$
$$+
\frac{1}{n!} \Bigl((\alpha-1)^n+\alpha^n \Bigr)=
\frac{1}{n!}\sum_{i=0}^{n} {n\choose i}^2 (\alpha-1)^{n-i}
\alpha^i=$$
$$=
\frac{1}{2^n n!}\sum_{i=0}^{n} {n\choose i}^2 (\gamma-1)^{n-i}
(\gamma+1)^i.$$
We took into account that $$\displaystyle { {n-1}\choose i } + { {n-1}\choose {i-1} } =
{n\choose i}$$ completing the proof of the left-hand equality in \eqref{th_2_1_1}.

The right-hand equality in \eqref{th_2_1_1}
follows from the identity
$$\sum_{i=0}^n {n\choose i}^2 t^i=(1-t)^n \chi_n\left(\frac{1+t}{1-t}\right) $$
(see \cite{prudnikov_1981}).
Let us put $t=\frac{\gamma - 1}{\gamma+1},$ then
$$(1-t)^n=2^n (\gamma+1)^{-n}, \quad \frac{1+t}{1-t}=\gamma.$$
Hence,
$$\mes(G)=
\frac{1}{2^n n!} \sum_{i=0}^{n} {n \choose i}^2 (\gamma - 1)^{n-i}
(\gamma + 1)^i =
$$
$$=
\frac{1}{2^n n!} \sum_{i=0}^{n} {n \choose i}^2 (\gamma+1)^{n-i}
(\gamma -1)^i =$$
$$
=\frac{1}{2^n n!} (\gamma+1)^n
\sum_{i=0}^{n} {n \choose i}^2 \Bigl( \frac{\gamma-1}
{\gamma +1}\Bigr)^i =
\frac{\chi_n(\gamma)}{n!}.$$

The  proof of the theorem is complete.
\hfill$\Box$

\msk
Let us give some simple examples. The set
$$E_{1,2}=\{ x\in {\mathbb R}: |x|+|1-x|\leq 2\}$$ is the segment
$\left[-\frac{1}{2},\frac{3}{2}\right]$ with the length
$${\rm mes_1}(E_{1,2}) =\frac{\chi_1(2)}{1!}=2.$$

\par
\begin{figure}[h]
\center{\includegraphics[scale=0.42]{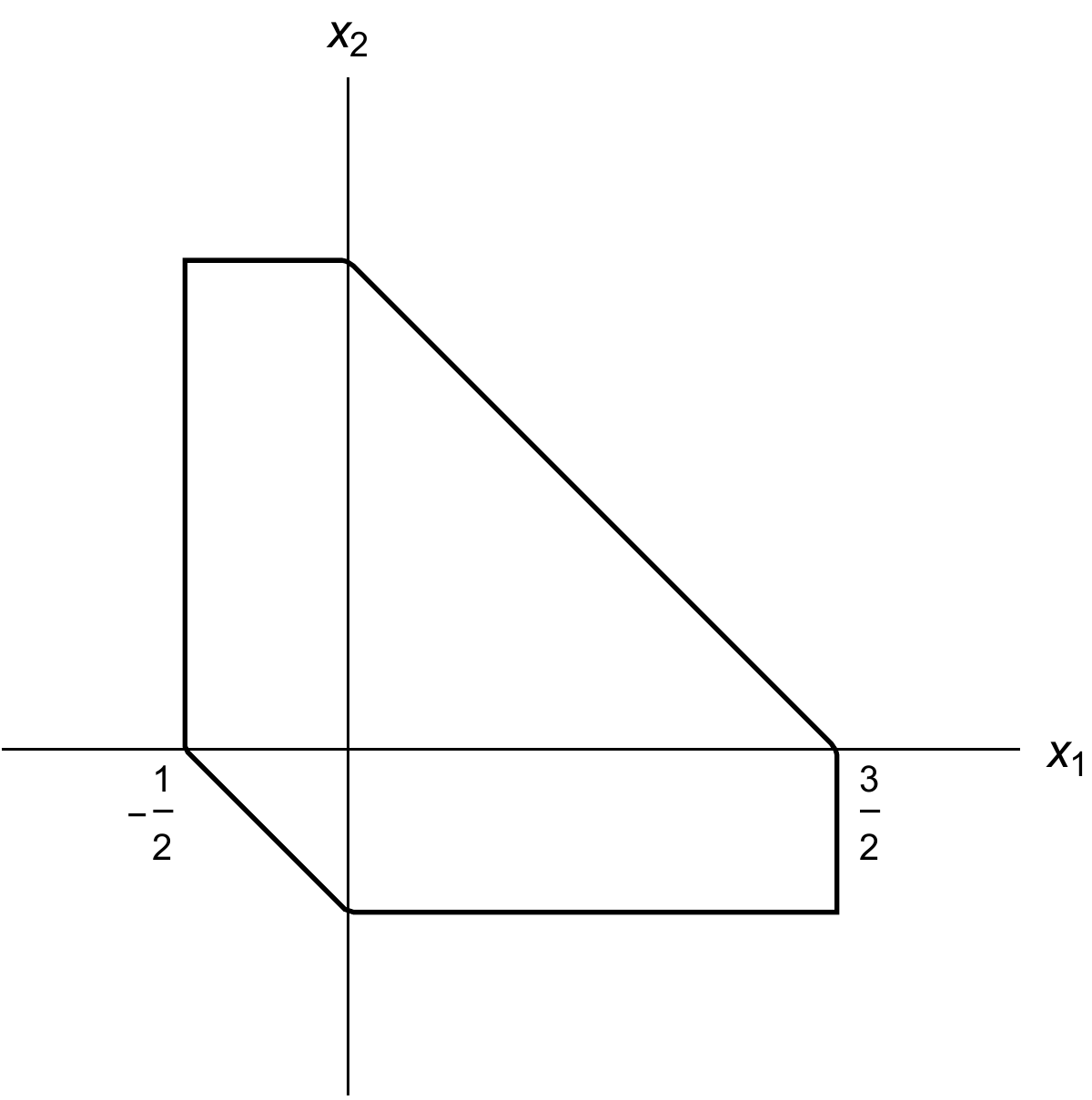}}
\caption{The set $E_{2,2}$}
\label{fig:nev_ukl_example_E_2_2}
\end{figure}

 The two-dimensional set
$
E_{2,2}=\{ x\in {\mathbb R}^2: |x_1|+|x_2|+|1-x_1-x_2|\leq	 2\}
$
is shown in Fig.~\ref{fig:nev_ukl_example_E_2_2}.

\par The area of this set is equal to
$${\rm mes_2}(E_{2,2}) =\frac{\chi_2(2)}{2!}=\frac{11}{4}.$$

\par
\begin{figure}[h!]
\center{\includegraphics[scale=0.48]{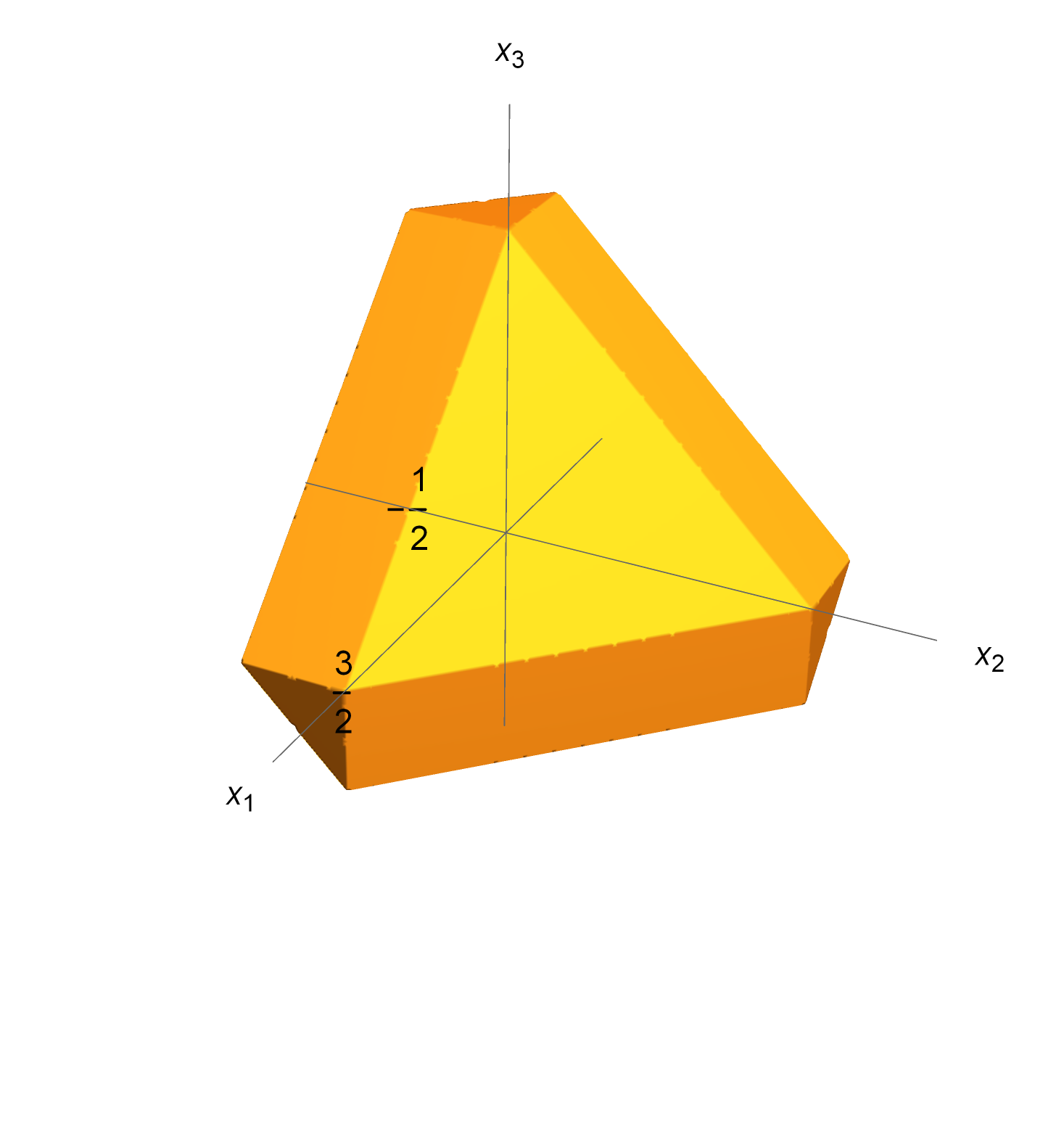}}
\vspace{-3cm}
\caption{The set $E_{3,2}$}
\label{fig:nev_ukl_example_E_3_2_one}
\end{figure}

 The three-dimensional domain
$$
E_{3,2}=\{ x\in {\mathbb R}^3: |x_1|+|x_2|+|x_3|+|1-x_1-x_2-x_3|\leq	 2\}
$$
is shown in Fig.~\ref{fig:nev_ukl_example_E_3_2_one}. It has the volume
$${\rm mes_3}(E_{3,2})=\frac{\chi_3(2)}{3!}=\frac{17}{6}.$$

There is an interesting open problem  related to the above-mentioned equality \eqref{th_2_1_1}. Essentially, 
this~equality
along with Rodrigues' formula and other well-known relations also gives a characterization of Legendre polynomials. This characterization is written via the volumes of some convex polyhedra. Namely, for $t\geq 1$
  \begin{equation}\label{Legendre}
\chi_n(t)=
n!\,{\rm mes}_n(E_{n,t}),
\end{equation}
where
$$E_{n,t}= \left\{ x\in{\mathbb R}^n :
\sum_{j=1}^n\left |x_j\right| +
\left|1- \sum_{j=1}^n x_j\right| \leq t \right\}.
$$
\smsk

The question arises about analogues of \eqref{Legendre} for other classes of orthogonal polynomials, such as
Chebyshev polynomials or, more generally, Jacobi polynomials: {\it Is the equality \eqref{Legendre} a particular case of a more general pattern?}

The author would be grateful for any information on this matter.

Let ua note that, from  \eqref{Legendre} and \eqref{reccurent_Legendre}, we have
$$
{\rm mes}_{n+1}(E_{n+1,t})=\frac{2n+1}{(n+1)^2}\,t\, {\rm mes}_n(E_{n,t})-\frac{1}{(n+1)^2}{\rm mes}_{n-1}(E_{n-1,t}).
$$

The direct establishing this recurrence relation for the measures of the set $E_{n,t}$ could provide a proof of
Theorem 5.1 different from the above.



\SECT{6. Inequalities $\theta_n(Q_n)>c\sqrt{n}$ and $\theta_n(B_n)>c\sqrt{n}$}{6}
\label{th_geq_c_sqrt_n_ineqs}
\addtocontents{toc}{6. Inequalities $\theta_n(Q_n)>c\sqrt{n}$ and $\theta_n(B_n)>c\sqrt{n}$
\hfill\thepage\par\VST}

\indent\par Based on Theorem 5.1, in this section we obtain lower bounds for the norms of the projector associated with linear interpolation on the cube and the ball.

First let us consider linear interpolation on the unit cube $Q_n=[0,1]^n$.

\smallskip
{\bf Theorem 6.1.} {\it Suppose $P: C(Q_n)\to \Pi_1\left({\mathbb R}^n\right)$ is an arbitrary interpolation
projector. Then for~the~corresponding simplex $S\subset Q_n$ and the node matrix
${\bf A}$ we have
\begin{equation}\label{th_2_3_formulation}
\|P\|_{Q_n}
\geq
\chi_n^{-1} \left(\frac{n!}{|\det({\bf A})}\right)=\chi_n^{-1} \left(\frac{1}{\vo(S)}\right).
\end{equation}
}

\smallskip
{\it Proof.}
For each $i=1,$ $\ldots,$ $n$, let us  subtract from the
$i$th row of ${\bf A}$ its $(n+1)$th row.
Denote by ${\bf B}$ the submatrix of order $n$ which stands in the first
$n$ rows and columns of the result matrix. Then
$$
|\det({\bf B})|=|\det({\bf A})|= n!\vo(S)\leq n!\nu_n
$$
so that
\begin{equation}\label{proof_111}
\frac{|\det({\bf B})|}{n!\nu_n}\leq 1.
\end{equation}

Let $x^{(j)}$ be the vertices and $\lambda_j$ be the basic Lagrange polynomials of  simplex $S$.
Since $\lambda_1(x),$ $\ldots$ , $\lambda_{n+1}(x)$ are the barycentric coordinates of a point $x$, we have
$$
\|P\|_{Q_n}=\max_{x\in Q_n} \sum_{j=1}^{n+1}|\lambda_j(x)|
=\max
\left\{ \sum_{j=1}^{n+1}|\beta_j|:
 \, \sum_{j=1}^{n+1}\beta_j=1,
\,
\sum_{j=1}^{n+1}\beta_jx^{(j)}\in Q_n
\right\}.
$$

Let us replace $\beta_{n+1}$ with the equal value
$1-\sum\limits_{j=1}^n
\beta_j.$ The condition
$\sum\limits_{j=1}^{n+1}\beta_jx^{(j)}$ $\in Q_n$ is equivalent to
$$
\sum\limits_{j=1}^{n}\beta_j(x^{(j)}-x^{(n+1)})\in Q^\prime:=Q_n-x^{(n+1)}.
$$
Hence,
\begin{equation}\label{norm_P_barycentric_from1_ton}
\|P\|_{Q_n}=\max\left \{ \sum_{j=1}^{n} \left|\beta_j\right| +
\left|1-\sum_{j=1}^n\beta_j\right|
\right \}
\end{equation}
where the maximum is taken over all
$\beta_j$ such that
$\sum\limits_{j=1}^n\beta_j(x^{(j)}-x^{(n+1)})
\in Q^\prime.$ Clearly,  
$$
\vo(Q^\prime)=\vo(Q_n)=1.
$$

Let us consider the nondegenerate linear operator $F:{\mathbb R}^n \to {\mathbb R}^n,$
which maps the point $\beta=(\beta_1,\ldots,\beta_n)$ to  the point
$x=F(\beta)$
according to the rule $$x=\sum_{j=1}^n \beta_j \left(x^{(j)}-x^{(n+1)}\right).$$
We have the matrix equality
$F(\beta)=(\beta_1,\ldots,\beta_n) {\bf B},$
where $\bf B$ is the above introduced matrix of order $n$ with elements $b_{ij}=x_j^{(i)}-x_j^{(n+1)}.$
Put $$\gamma^*:=\chi_n^{-1}\left(\displaystyle \frac{n!}{|\det {\bf B}|}\right).$$
This definition is correct since by \eqref{proof_111}
$$\displaystyle\frac{n!}{|\det({\bf B})|}\geq \frac{1}{\nu_n}\geq 1.$$
Note that
$$
\chi_n(\gamma^*)=\frac{n!}{|\det({\bf B})|}.
$$

Given $\gamma\geq 1$, let us introduce a set
$$
E_{n,\gamma}:= \left\{ \beta=(\beta_1,\ldots,\beta_n)\in{\mathbb R}^n :
\sum_{j=1}^n |\beta_j| +
\left|1- \sum_{j=1}^n \beta_j\right|  \leq \gamma \right\}.
$$

Let us show that $Q^\prime \not \subset F(E_{n,\gamma})$ provided $\gamma < \gamma^*.$ Indeed, this is immediate from Theorem 5.1 because, thanks to this theorem,
$$
\mes(F(E_{n,\gamma})) <
\mes(F(E_{n,\gamma^*}))=|\det {\bf B}|\cdot \mes(E_{n,\gamma^*})=
$$
$$=
|\det {\bf B}|\cdot
\frac{\chi_n(\gamma^*)}{n!}=1=\mes(Q^\prime).
$$

Thus, for every $\varepsilon > 0$ there exists a point
$x^{(\varepsilon)}$ with the properties:
$$x^{(\varepsilon)}
\sum\beta_j^{(\varepsilon)} \left(x^{(j)}-x^{(n+1)}\right)\in Q^\prime, \quad
\left|\sum\beta_j^{(\varepsilon)}\right|
+\left|1-\sum\beta_j^{(\varepsilon)}\right|
\geq \gamma^*-\varepsilon.$$
In view of \eqref{norm_P_barycentric_from1_ton} this gives
$\|P\|_{Q_n}\geq \gamma^*-\varepsilon.$
Since~$\varepsilon>0$ is an arbitrary, we obtain
$$\|P\|_{Q_n} \geq \gamma^*=
\chi_n^{-1} \left(\frac{n!}{|\det({\bf B})|}\right)=
\chi_n^{-1} \left(\frac{n!}{|\det({\bf A})|}\right)=\chi_n^{-1} \left(\frac{1}{\vo(S)}\right).$$
The theorem is proved.
\hfill$\Box$

\msk
{\bf Corollary 6.1.} {\it For each $n$,
\begin{equation}\label{corol_2_1}
\theta_n(Q_n)\geq
\chi_n^{-1} \left(\displaystyle\frac{1}{\nu_n}\right).
\end{equation}}

\smallskip
The Stirling formula (see, e.\,g., \cite{fiht_2001})
$$
n!=\sqrt{2\pi n}\left(\frac{n}{e}\right)^n
e^{\frac{\zeta_n}{12n}},~~~
0<\zeta_n<1,
$$
yields
\begin{equation}\label{n_fact_ineqs}
\sqrt{2\pi n}\left(\frac{n}{e}\right)^n<n!<\sqrt{2\pi n}\left(\frac{n}{e}\right)^n
e^{\frac{1}{12n}}.
\end{equation}

Let us see that for every $n>1$ the following inequality
\begin{equation}\label{chi_ineq}
\chi_n^{-1}(s)>\left( \frac{s}{ {n\choose{\lfloor n/2\rfloor}} }\right)^{1/n}
\end{equation}
holds. In fact, if $t\geq 1$ and $n>1$, then, according to \eqref{th_2_1_1},
$$\chi_n(t)=
\frac{1}{2^n}\sum_{i=0}^{n} {n\choose i}^2 (t-1)^{n-i}(t+1)^i<$$
$$
<{n\choose{\lfloor n/2 \rfloor}} \cdot \frac{1}{2^n}\sum_{i=0}^{n}
{n\choose i} (t-1)^{n-i}(t+1)^i=
$$
$$= {n\choose{\lfloor n/2 \rfloor}}\cdot (2t)^n \cdot 2^{-n}=
{n\choose{\lfloor n/2 \rfloor} }\, t^n
$$
proving \eqref{chi_ineq}.

For even $n$,
$\displaystyle{n\choose{\lfloor n/2 \rfloor} }=
{n\choose{n/2}}=(n!)/((n/2)!)^2,$ therefore,
\begin{equation}\label{chi_ineq_2}
\chi_n^{-1}(s)>\left( \frac{s\left((n/2)!\right)^2}
{ n!}\right)^{1/n}.
\end{equation}
If $n$ is odd, then
$${n\choose{\lfloor n/2 \rfloor}}
=\frac{n!}{  \frac{n+1}{2} ! \frac{n-1}{2}!  }$$
and \eqref{chi_ineq}
has the form
\begin{equation}\label{chi_ineq_3}
\chi_n^{-1}(s)>\left( \frac{s\frac{n+1}{2}!\frac{n-1}{2}!}
{ n!}\right)^{1/n}.
\end{equation}

\smallskip
{\bf Theorem 6.2.}
{\it For all $n\in {\mathbb N}$,
\begin{equation}\label{th_2_4_form}
\theta_n(Q_n)>\frac{\sqrt{n-1}}{e}.
\end{equation}
}

\smallskip
{\it  Proof.} The case $n=1$ is trivial.  If $n>1,$
let us  apply \eqref{chi_ineq}--\eqref{chi_ineq_3}
and the Hadamard inequality (see~\eqref{adamar_clements_lindstrem})
$$\nu_n\leq\frac{\left(n+1\right)^{(n+1)/2}}{2^nn!}.$$
For even $n$, from \eqref{corol_2_1},
the left-hand inequality in \eqref{n_fact_ineqs}, and
  \eqref{chi_ineq_2},
we get
$$\theta_n(Q_n)\geq \chi_n^{-1}\left(\frac{1}{\nu_n}\right)\geq
\chi_n^{-1}\left( \frac{2^nn!}{(n+1)^{(n+1)/2}}\right)>
2\left(\frac{[(n/2)!]^2}{(n+1)^{(n+1)/2}}\right)^{1/n}>$$
$$>
\frac{2}{(n+1)^{1/2+1/(2n)}}\left( \sqrt{\pi n} \left(\frac{n}{2e}\right)^{n/2}
\right)^{2/n}=
\frac{ \left(\pi n\right)^{1/n}n }{e(n+1)^{1/2+1/(2n)}}>
\frac{\sqrt{n-1}}{e}.$$
If $n$ is odd, then, thanks to \eqref{chi_ineq_3},
$$\theta_n(Q_n)\geq \chi_n^{-1}\left(\frac{1}{\nu_n}\right)\geq
\chi_n^{-1}\left( \frac{2^nn!}{(n+1)^{(n+1)/2}}\right)>
\left( \frac{2^n \frac{n+1}{2}!\frac{n-1}{2}! } {(n+1)^{(n+1)/2}}
\right)^{1/n}>$$
$$>
2\left( \frac{ \pi\sqrt{n^2-1}\left(n^2-1\right)^{(n-1)/2}
(n+1) }{(2e)^n} \right)^{1/n}=$$
$$=
\frac{1}{e}\pi^{1/n}\sqrt{n-1}(n+1)^{1/(2n)}>
\frac{\sqrt{n-1}}{e}$$
proving that \eqref{th_2_4_form} holds for all $n$.
\hfill$\Box$

\smallskip
In some situations, the estimates of Theorem 2.4 can be improved \cite{nevskii_monograph}. In particular, we prove 
that $\theta_n>\sqrt{n}/{e}$ if $n$ is even, and
$\displaystyle \theta_n>n/(e\sqrt{n-1})$ provided $n>1$ and $n\equiv 1({\rm mod}~4)$.

Let us also note that, thanks to \eqref{th_2_4_form}, the inequality $\theta_n(Q_n)> c\sqrt{n}$ holds with some $c$, $0<c<1$.
A~suitable estimate is
\begin{equation}\label{theta_n_lt_sqrt_n}
\theta_n(Q_n)>\frac{2\sqrt{2}}{3e} \, \sqrt{n}.
\end{equation}

Indeed, if $n\leq 8$, then the right-hand side of \eqref{theta_n_lt_sqrt_n} is less than 1,  while for
$n\geq 9$
$$ \theta_n(Q_n)>\frac{\sqrt{n-1}}{e}\geq \frac{2\sqrt{2}}{3e} \, \sqrt{n}.$$
Notice that
$\displaystyle\frac{2\sqrt{2}}{3e} = 0.3468 \ldots$
\bsk

Later this approach was extended to linear interpolation on the unit ball $B_n$, see~\cite{nevskii_mais_2019_26_3}.  Let us present some results obtained in this direction.

A regular simplex inscribed into an $n$-dimensional ball has the maximal volume among all simplices contained in this ball. Furthermore, there are no other simplices with this property  (see \cite{fejes_tot_1964}, \cite{slepian_1969}, and
\cite{vandev_1992}).

Let $\varkappa_n=\vo(B_n)$, and let $\sigma_n $ be the volume of a regular simplex inscribed into $B_n$.

\smallskip
{\bf Theorem 6.3.}
{\it Let
$P:C(B_n)\to\Pi_1({\mathbb R}^n)$ be an interpolation projector. Then, for the corresponding simplex
$S\subset B_n$ and the node matrix ${\bf A}$, we have
\begin{equation}\label{norm_P_vol_S_ineq}
\|P\|_{B_n}
\geq
\chi_n^{-1}
\left(\frac{n!\varkappa_n}{|\det({\bf A})|}\right)=
\chi_n^{-1}
\left(\frac{\varkappa_n}{\vo(S)}\right).
\end{equation}
}

\smallskip


\smallskip
{\bf Corollary 6.2.}
{\it For each $n$,
\begin{equation}\label{theta_n_chi_n_kappa_sigma_ineq_1}
\theta_n(B_n)
\geq
\chi_n^{-1}
\left(\frac{\varkappa_n}{\sigma_n}\right).
\end{equation}
}

\smallskip
{\it Proof.}
Let $P: C(B_n) \to \Pi_1({\mathbb R}^n)$ be an arbitrary interpolation projector. The
corresponding simplex $S\subset B_n$ satisfies  $\vo(S)\leq\sigma_n$. Then, thanks to inequality (\ref{norm_P_vol_S_ineq}),
$$
\|P\|_{B_n}
\geq
\chi_n^{-1}
\left(\frac{\varkappa_n}{\vo(S)}\right)
\geq
\chi_n^{-1}
\left(\frac{\varkappa_n}{\sigma_n}\right)
$$
proving (\ref{theta_n_chi_n_kappa_sigma_ineq_1}).
\hfill$\Box$

\smallskip
It is known that
\begin{equation}\label{kappa_n_sigma_n_whole_1}
\varkappa_n=\frac{\pi^
{\frac{n}{2}}}
{\Gamma\left(\frac{n}{2}+1\right)},\qquad
\sigma_n=\frac{1}{n!}\sqrt{n+1}\left(\frac{n+1}{n}\right)^{\frac{n}{2}},
\end{equation}
\begin{equation}\label{kappa_n_even_and_odd_1}
\varkappa_{2k}=\frac{\pi^{k}}{k!},\qquad
\varkappa_{2k+1}=\frac{2^{k+1}\pi^{k}}{(2k+1)!!}=
\frac{2(k!)(4\pi)^k}{(2k+1)!}.
\end{equation}
These equalities and inequality (\ref{theta_n_chi_n_kappa_sigma_ineq_1}) imply the following low bound for the constant $\theta_n(B_n)$. 

\smallskip
{\bf Corollary 6.3.}
{\it For every $n$,
\begin{equation}\label{theta_n_chi_n_through_n_ineq}
\theta_n(B_n)
\geq
\chi_n^{-1}
\left(\frac{\pi^{\frac{n}{2}}n!}{\Gamma\left(\frac{n}{2}+1\right)\sqrt{n+1}
\left(\frac{n+1}{n}\right)^{\frac{n}{2}}}
\right).
\end{equation}
If $n=2k$, then  (\ref{theta_n_chi_n_through_n_ineq}) is equivalent to the inequality
\begin{equation}\label{theta_n_chi_n_through_n_is_2k_ineq}
\theta_{2k}(B_{2k})
\geq
\chi_{2k}^{-1}
\left(\frac{\pi^{k}(2k)!}{k!\sqrt{2k+1}
\left(\frac{2k+1}{2k}\right)^k}
\right).
\end{equation}
For $n=2k+1$ we have
\begin{equation}\label{theta_n_chi_n_through_n_is_2k_plus_1_ineq}
\theta_{2k+1}(B_{2k+1})
\geq
\chi_{2k+1}^{-1}
\left(\frac{2(k!)(4\pi)^{k}}{\sqrt{2k+2}
\left(\frac{2k+2}{2k+1}\right)^{\frac{2k+1}{2}}}
\right).
\end{equation}
}

\smallskip
{\it Proof.} It is sufficient to apply (\ref{theta_n_chi_n_kappa_sigma_ineq_1}),
(\ref{kappa_n_sigma_n_whole_1}), and
(\ref{kappa_n_even_and_odd_1}).
\hfill$\Box$
\msk

Also we will need the following estimates
which were proved in \cite[Section\,3.4.2]{nevskii_monograph}:
\begin{equation}\label{chi_2k_2k_plus_1_ineqs}
\chi_{2k}^{-1}(s)>\left( \frac{(k!)^2s}
{ (2k)!}\right)^{\frac{1}{2k}}, \qquad
\chi_{2k+1}^{-1}(s)>\left( \frac{(k+1)!k!s}
{ (2k+1)!}\right)^{\frac{1}{2k+1}}.
\end{equation}

\smallskip
{\bf Theorem 6.4.}
{\it
 There exists an absolute constant $c>0$ such that
\begin{equation}\label{theta_n_B_n_gt_c_sqrt_n_1}
\theta_n(B_n)>c\sqrt{n}.
\end{equation}
}

\smallskip
In \cite{nevskii_mais_2021_28_2}  we have shown that
the inequality
(\ref{theta_n_B_n_gt_c_sqrt_n_1}) takes place  with the constant

\begin{equation}\label{C1}
c=
\frac{  \sqrt[3]{\pi} }{\sqrt{12e}\cdot\sqrt[6]{3} } =  0.2135...\,.
\end{equation}

\smallskip
{\bf Corollary 6.4.}
{\it   $\theta_n(B_n)\asymp \sqrt{n}$.
}

\smallskip
{\it Proof.}
By the results of \cite{nev_ukh_mais_2019_26_2},  $\theta_n(B_n)\leq \sqrt{n+1}.$
Consequently, the lower estimate $\theta_n(B_n)>c\sqrt{n}$
is~precise with respect to dimension $n$.
\hfill$\Box$

\msk
{\bf Corollary 6.5.}
{\it  Let $P:C(B_n)\to\Pi_1({\mathbb R}^n)$ be the
interpolation projector whose nodes coincide with vertices of a regular simplex being inscribed into the boundary sphere $\|x\|=1$.
Then $\|P\|_{B_n}\asymp \theta_n(B_n)$.}

\smallskip
{\it Proof.}
Since $\sqrt{n}\leq \|P\|_{B_n}\ \leq \sqrt{n+1}$,   it remains to
utilize the previous corollary. \hfill$\Box$

\bsk
Until 2023, the equality $\|P\|_{B_n}=\theta_n(B_n)$, where $P$ is the projector from Corollary 6.5, remained proven
only for $1\leq n\leq 4$. As it turned out later, this equality is true for every $n$. Moreover, if~$n>1$, then, in addition to the inequality (\ref{theta_n_B_n_gt_c_sqrt_n_1}) with $c$ defined by (\ref{C1}), we proved that 
$\theta_{n}(B_{n})>\sqrt{n}$. 
However, the latter only became clear after finding the exact value of $\theta_n(B_n)$. See Section 9.

\msk
The approach proposed above can be also applied to the case of an arbitrary convex body $K\subset {\mathbb R}^n.$

\smallskip
{\bf Theorem 6.5.}
{\it  For an arbitrary interpolation projector
$P:C(K)\to\Pi_1({\mathbb R}^n)$,
\begin{equation}\label{norm_P_vol_S_for_K_ineq}
\|P\|_{K}
\geq
\chi_n^{-1}
\left(\frac{n!\vo(K)}{|\det({\bf A})|}\right)=
\chi_n^{-1}
\left(\frac{\vo(K)}{\vo(S)}\right).
\end{equation}
Here $S\subset K$ is the corresponding simplex and  ${\bf A}$ is the corresponding node matrix.
}

\msk
We let $\simp(K)$ denote the maximum volume of a nondegenerate simplex $S$ with vertices in $K$. Thus, in these settings, $\nu_n=\simp(Q_n)$ and $\sigma_n=\simp(B_n)$.

\smallskip
{\bf Corollary 6.6.}
{\it Let $K$ be an arbitrary convex body in ${\mathbb R}^n$. Then
\begin{equation}\label{theta_n_K_ineq_thru_xi_minus_one}
\theta_n(K)
\geq
\chi_n^{-1}
\left(\frac{\vo(K)}{\simp(K)}\right).
\end{equation}
}

\smallskip
This inequality is immediate from \eqref{norm_P_vol_S_for_K_ineq}.

\smsk
For a simplex $S$ of the maximum volume in $K$, we have $K\subset (n+2)S,$ see
\eqref{S_sub_K_sub_n_plus_2_S}. Therefore, the~ratio $\vo(K)/\simp(K)$ in
\eqref{theta_n_K_ineq_thru_xi_minus_one} is bounded by $(n+2)^n.$

Note that the approach given in this section works also for an arbitrary (not necessarily convex) compact.



\SECT{7. The norm $\|P\|_{B}$ for an inscribed regular simplex}{7}
\label{norm_P_for_insc_reg_simplex}
\addtocontents{toc}{7. The norm $\|P\|_{B}$ for an inscribed regular simplex\hfill\thepage\par\VST}

\indent\par Let $B=B(x^{(0)};R)$ and let $P:C(B)\to \Pi_1({\mathbb R}^n)$ be
the~interpolation projector with the nodes  $x^{(j)}\in B$, $j=1,...,n+1$.
Let $S$ be the simplex with vertices  $x^{(j)}$ and let
$\lambda_j(x)=l_{1j}x_1+ \ldots+l_{nj}x_n+l_{n+1,j}$ be
the basic Lagrange polynomials of $S$. In this case
$$\|P\|_{B}=
\max_{x\in B}\sum_{j=1}^{n+1}
|\lambda_j(x)|
$$
(see (\ref{norm_P_intro_cepochka}) for $K=B$).
In \cite{nev_ukh_mais_2019_26_2}   another formula for $\|P\|_{B}$ was obtained:

$$\|P\|_B=
\max\limits_{f_j=\pm 1} \left[
R \left(\sum_{i=1}^n\left(\sum_{j=1}^{n+1} f_jl_{ij}\right)^2\right)^{1/2}
+\left|\sum_{j=1}^{n+1}f_j\left(\sum_{i=1}^n l_{ij}x_i^{(0)}+l_{n+1,j}\right)\right|
\right]=$$
\begin{equation}\label{norm_P_B_formula}
=
\max\limits_{f_j=\pm 1} \left[
R \left(\sum_{i=1}^n\left(\sum_{j=1}^{n+1} f_jl_{ij}\right)^2\right)^{1/2}
+\left|\sum_{j=1}^{n+1}f_j\lambda_j(x^{(0)})\right|
\right].
\end{equation}

\noindent If $c(S)=x^{(0)}$, then $\lambda_j(x^{(0)})=1/(n+1)$ which implies the following
simpler formula for $\|P\|_B$:
\begin{equation}\label{norm_P_B_formula_c_eq_x0}
\|P\|_B=
\max\limits_{f_j=\pm 1} \left[
R \left(\sum_{i=1}^n\left(\sum_{j=1}^{n+1} f_jl_{ij}\right)^2\right)^{1/2}
+\frac{1}{n+1}\left|\sum_{j=1}^{n+1}f_j\right|
\right].
\end{equation}

Now, suppose that $S$ is a regular simplex inscribed into the $n$-dimensional ball
$B=B(x^{(0)};R)$ and
$P:C(B)\to \Pi_1({\mathbb R}^n)$ is the corresponding  interpolation projector.
Clearly, $\|P\|_B$ does not depend on the center $x^{(0)}$ and
the radius $R$ of the ball and  on the choice of a regular simplex
inscribed into that ball.
In other words, $\|P\|_B$ is a function of only  dimension $n$. The exact expression for $\|P\|_B$ which we present below was found in \cite{nev_ukh_mais_2019_26_2}.

We let $\psi:[0,n+1]\to\R$ denote a function defined by
\begin{equation}\label{psi_function_modulus}
\psi(t):=\frac{2\sqrt{n}}{n+1}\Bigl(t(n+1-t)\Bigr)^{1/2}+
\left|1-\frac{2t}{n+1}\right|.
\end{equation}

Let
\begin{equation}\label{a_formula}
a=a_n=\left\lfloor\frac{n+1}{2}-\frac{\sqrt{n+1}}{2}\right\rfloor.
\end{equation}

{\bf Theorem 7.1.} {\it The following relations hold:
\begin{equation}\label{norm_P_reg_formula}
\|P\|_B=\max\{\psi(a),\psi(a+1)\},
\end{equation}
\begin{equation}\label{norm_P_reg_ineqs}
\sqrt{n}\leq \|P\|_B \leq \sqrt{n+1}.
\end{equation}
The equality $\|P\|_B=\sqrt{n}$ is true only for $n=1$.
The equality $\|P\|_B=\sqrt{n+1}$ holds if and only if
$\sqrt{n+1}$ is an integer.
}

{\it Proof.} First  we  prove (\ref{norm_P_reg_formula}).
If $n=1$, then
$\psi(t)=\sqrt{t(2-t)}+|1-t|$, $a=0$, $\psi(a)=\psi(a+1)=1.$
Since $\|P\|_B=1$, the equality (\ref{norm_P_reg_formula}) is true.

Let $n\geq 2$.
Consider the simplex $S$
with vertices
$$x^{(1)}=e_1, \ \  \ldots, \ \ x^{(n)}=e_n, \quad
x^{(n+1)}=
\left(\frac{1-\sqrt{n+1}}{n},\ldots,\frac{1-\sqrt{n+1}}{n}\right).
$$
Since  the length of any edge of $S$ is equal to
$\sqrt{2}$,  this simplex is regular.
The simplex is inscribed into the ball
$B=B(x^{(0)};R)$, where
$$x^{(0)}=c(S)=\left(\frac{1-\sqrt{\frac{1}{n+1}}}{n},\ldots,
\frac{1-\sqrt{\frac{1}{n+1}}}{n}\right),
\quad R=\sqrt{\frac{n}{n+1}}.$$
Note that the $(n+1)$th vertex of $S$ is obtained by shifting the zero-vertex of the simplex  $x_i\geq 0$, $\sum x_i\leq 1$ in direction from the hyperplane $\sum x_i=1$.
It is important that $S$ is invariant with respect to changing the order of coordinates.
It is enough to find $\|P\|_B$ for this simplex.

The corresponding matrices ${\bf A}$ and ${\bf A}^{-1}$ are
\begin{equation}\label{A_A_minus_1_for_regular_S}
{\bf A}=
\left( \begin{array}{ccccc}
1&0&\ldots&0&1\\
0&1&\ldots&0&1\\
\vdots&\vdots&\vdots&\vdots&\vdots\\
0&0&\ldots&1&1\\
-\tau&-\tau&\ldots&-\tau&1\\
\end{array}
\right), \quad
{\bf A}^{-1} =
\frac{1}{\sqrt{n+1}}
\left( \begin{array}{ccccc}
\sigma&-\tau&\ldots&-\tau&-1\\
-\tau&\sigma&\ldots&-\tau&-1\\
\vdots&\vdots&\vdots&\vdots&\vdots\\
-\tau&-\tau&\ldots&\sigma&-1\\
\tau&\tau&\ldots&\tau&1\\
\end{array}
\right).
\end{equation}
Here
\begin{equation}\label{sigma_tau_formulae}
\sigma:=\frac{(n-1)\sqrt{n+1}+1}{n}, \quad \tau:=\frac{\sqrt{n+1}-1}{n}.
\end{equation}
Since $c(S)=x^{(0)}$, for calculation $\|P\|_B$ we can use
(\ref{norm_P_B_formula_c_eq_x0}) with
 $R=\sqrt{n/(n+1)}$:
$$
\|P\|_B=
\max\limits_{f_j=\pm 1} \left[
\sqrt{\frac{n}{n+1}}
 \left(\sum_{i=1}^n\left(\sum_{j=1}^{n+1} f_jl_{ij}\right)^2\right)^{1/2}
+\frac{1}{n+1}\left|\sum_{j=1}^{n+1}f_j\right|
\right].
$$
Here $l_{ij}$ are the elements of ${\bf A}^{-1}$ (see (\ref{A_A_minus_1_for_regular_S})).

Let $k$ be the number of $f_j$ equal to $-1$.
Then the number of $f_j$ equal to $1$ is $n+1-k$.
The simplex $S$ does not change with renumerating of coordinates, so,
we may assume that
$f_1=\ldots=f_k=-1,$ $f_{k+1}=\ldots=f_{n+1}=1$.
Since the function being maximized does not change when the signs of
$f_j$ change simultaneously, we can consider only the interval $1 \leq k \leq \frac{n + 1}{2}$.
Thus,
\begin{equation}\label{P_norm_f_j_long}
\|P\|_B=
\max\limits_{1\leq k\leq \frac{n+1}{2}} \left[
\sqrt{\frac{n}{n+1}}
 \left(\sum_{i=1}^n\left(\sum_{j=1}^{n+1} f_jl_{ij}\right)^2\right)^{1/2}
+\frac{n+1-2k}{n+1}
\right],
\end{equation}
where $f_j=-1$ for $1\leq j\leq k$ and $f_j=1$ for all other $j$.
The number $n+1-2k$ is equal to the sum under the sign of absolute value.
Taking into account the multiplier
$1/\sqrt{n+1}$
 in the equality for
${\bf A}^{-1}$, let us rewrite the value
$$
W=(n+1)\sum_{i=1}^{n}
\left(\sum_{j=1}^{n+1}f_jl_{ij}\right)^2.$$
Utilizing the explicit expressions for $l_{ij}$ we present this sum in the form
$$W=
(n+1)\sum_{i=1}^{k}\left(\sum_{j=1}^{n+1}f_jl_{ij}\right)^2+
(n+1)\sum_{i=k+1}^{n}\left(\sum_{j=1}^{n+1}f_jl_{ij}\right)^2
=W_1+W_2.$$
From
(\ref{A_A_minus_1_for_regular_S}) and the described distribution of values $f_j$, we get
$$W_1=\sum_{i=1}^k (-\sigma+(k-1)\tau-(n-k)\tau-1)^2=k(2k\tau-\alpha)^2,$$
$$W_2=\sum_{i=k+1}^n (k\tau+\sigma-(n-1-k)\tau-1)^2=(n-k)(2k\tau+\beta)^2.$$
Here $\alpha=\sigma+(n+1)\tau+1$, $\beta=\sigma-(n-1)\tau-1$.
It follows from (\ref{sigma_tau_formulae}) that $\alpha=2\sqrt{n+1},$ $\beta=0$, therefore,
$$W=
4k(k\tau-\sqrt{n+1})^2+(n-k)\cdot4k^2\tau^2=k^2(-8\sqrt{n+1}\tau+4n\tau^2)+4k(n+1)=$$
$$=-4k^2+4k(n+1)=4k(n+1-k).$$
This gives
$$
\|P\|_B=
\max\limits_{1\leq k\leq \frac{n+1}{2}} \left[
\sqrt{\frac{n}{n+1}}
 \left(\frac{1}{n+1}W\right)^{1/2}
+\frac{n+1-2k}{n+1}\right]=$$
$$=\max\limits_{1\leq k\leq \frac{n+1}{2}} \left[
\frac{2\sqrt{n}}{n+1}
\bigl(k(n+1-k)\bigr)^{1/2}+1-\frac{2k}{n+1}\right].$$
Recalling (\ref{psi_function_modulus}), we obtain
\begin{equation}\label{norm_P_eq_max_psi}
\|P\|_B=
\max\limits_{1\leq k\leq \frac{n+1}{2}} \psi(k).
\end{equation}
It remains to show that the last maximum is equal to the largest of numbers
$\psi(a)$ and $\psi(a+1)$, where
$a=\left\lfloor\frac{n+1}{2}-\frac{\sqrt{n+1}}{2}\right\rfloor.$
To do this, we analyze the behavior of $\psi(t)$ over the whole interval $ [0, n + 1] $.

The function
$$\psi(t)=\frac{2\sqrt{n}}{n+1}\Bigl(t(n+1-t)\Bigr)^{1/2}+
\left|1-\frac{2t}{n+1}\right|, \quad 0\leq t\leq n+1,$$
has the following properties:
$$
\psi(t)>0,~~~
\psi(0)=\psi(n+1)=1,~~~
\psi\left(\frac{n+1}{2}\right)=\sqrt{n},~~~
\psi(n+1-t)=\psi(t).
$$
The graph of $\psi(t)$ is symmetric with respect to the straight line $t=\frac{n+1}{2}$.
On each half of $[0,n+1]$ the~function
$\psi(t)$ is concave as a sum of two concave functions.
Indeed, for $0\leq t\leq \frac{n+1}{2}$,
$$
\psi(t)=\varphi_1(t)+\varphi_2(t), \quad
\varphi_1(t)=\frac{2\sqrt{n}}{n+1}\Bigl(t(n+1-t)\Bigr)^{1/2}, \quad
\varphi_2(t)=1-\frac{2t}{n+1},
$$
where $\varphi_1(t)$ is concave as a superposition of the concave function
$t(n+1-t)$ and the increasing concave function $\sqrt{t}$, while
$\varphi_2(t)$ is a linear function.
The derivation $\psi^\prime(t)$ is equal to zero only in two points
$$t_-=\frac{n+1}{2}-\frac{\sqrt{n+1}}{2}, \quad
t_+=\frac{n+1}{2}+\frac{\sqrt{n+1}}{2}.$$

These points lie strictly inside the segments
$\left[0,\frac{n+1}{2}\right]$ and
$\left[\frac{n+1}{2},n+1\right]$ respectively.
From the conca\-vity of $\psi(t)$ on each of these segments,
it follows that
$$\max_{0\leq \psi(t)\leq n+1} \psi(t)=\psi(t_-)=\psi(t_+)=\sqrt{n+1}.$$
Moreover, $ t_-$ and $t _+ $ are the only maximum points on the left and right halves of $[0,n+1]$.
Hence, $\psi(t)$ increases for $0\leq t\leq t_-$ and decreases for $t_-\leq t\leq \frac{n+1}{2}$.
On the left half $\psi(t)$ behaves symmetrically: it~increases for  $\frac{n+1}{2}\leq t\leq t_+$
and decreases for $t_+\leq t\leq n+1.$

Let us now consider only the segment $\left[0,\frac{n+1}{2}\right]$.
Since $a=\left\lfloor t_-\right\rfloor$, then always $a\leq t_-<a+1$.
Let $k$ be a whole number, $1\leq k\leq \frac{n+1}{2}$.
If $k<a$, then $\psi(k)<\psi(a)$, while $\psi(k)<\psi(a+1)$ provided $k>a+1$.
Taking into account
(\ref{norm_P_eq_max_psi}),
we have
$$\|P\|_B=
\max\limits_{1\leq k\leq \frac{n+1}{2}} \psi(k)=\max \{\psi(a),\psi(a+1)\}.
$$
The equality (\ref{norm_P_reg_formula}) is proved.

We now turn to inequalities (\ref{norm_P_reg_ineqs}). We have already obtained the right one:
$$\|P\|_B=\max_{k\in {\mathbb Z}: \,1\leq k\leq\frac{n+1}{2}} \psi(k) \leq
\max_{1\leq t\leq\frac{n+1}{2}} \psi(t)=
\psi(t_-)=\sqrt{n+1}.$$

Let us describe such dimensions $n$ that $\|P\|_B=\sqrt{n+1}$.
These dimensions
are characterized by~the~fact that the non-strict inequality in the last relation
becomes an equality, i.e., the number $t_-$ is integer.

Note that both $t_-$ and $t_+$  are integer if and only if $\sqrt{n+1}$  is an integer.
For example, assume that 
$$
t_-=\frac{n+1}{2}-\frac{\sqrt{n+1}}{2}=d\in {\mathbb Z}.
$$
Then $\sqrt{n+1}=n+1-2d$ is an integer. On the contrary, if $\sqrt{n+1}=m\in{\mathbb Z}$,
then $m(m-1)$ is an even number and $t_-=\frac{m(m-1)}{2}$ is an integer.
The statement for  $t_+$ can be proved similarly. Also it can be deduced from the previous statement, since $t_+-t_-=\sqrt{n+1}$.

We have obtained that for any dimension of the form  $n=m^2-1$, $m$
 is an integer, and only in these situations, the number $t_-$
is integer.  Consequently,
$$\|P\|_B=\max_{k\in {\mathbb Z}: \,1\leq k\leq\frac{n+1}{2}} \psi(k)
=\max_{1\leq t\leq\frac{n+1}{2}} \psi(t)=
\psi(t_-)=\sqrt{n+1}=m.$$
These equalities are equivalent to (\ref{norm_P_reg_formula}),
since in the cases considered
$a=\left\lfloor t_-\right\rfloor=t_-$ and $\|P\|_B=\psi(a)$.
For all other $n$, holds $\|P\|<\sqrt{n+1}.$
Indeed, if $n \neq m^2-1$, then $a<t_-<a+1$  and the maximum
 in~(\ref{norm_P_eq_max_psi}) is reached either for $k=a$ or for $k=a+1$.
For any $n\ne m^2-1$ the norm of $P$, i.e., the~maximum $\psi(k)$ for integer
$k\in \left[1,\frac{n+1}{2}\right]$, is strictly less than maximum of $\psi(t)$ over all this interval.

It remains to show that always $\|P\|_B\geq \sqrt{n}$ and for $n>1$
this equality is a strict one.
If $n>3$, then $\sqrt{n+1}>2$, hence,
$$a+1=
\left\lfloor\frac{n+1}{2}-\frac{\sqrt{n+1}}{2}\right\rfloor+1\leq
\frac{n+1}{2}-\frac{\sqrt{n+1}}{2}+1<\frac{n+1}{2}.$$
From (\ref{norm_P_reg_formula}) and properties of $\psi(t)$, we have
$$\|P\|_B=\max \{\psi(a),\psi(a+1)\}\geq \psi(a+1)>\psi\left(\frac{n+1}{2}\right)=
\sqrt{n} \qquad(n>3).
$$
Also, $\|P\|_B>\sqrt{n}$ provided $n=2$ and $n=3$.
Therefore, $\|P\|=\sqrt{n}$ only for $n=1$, and the proof of~the~theorem is complete.
\hfill$\Box$

\medskip

\begin{figure}[h!]
 \centering
\includegraphics[width=11cm]{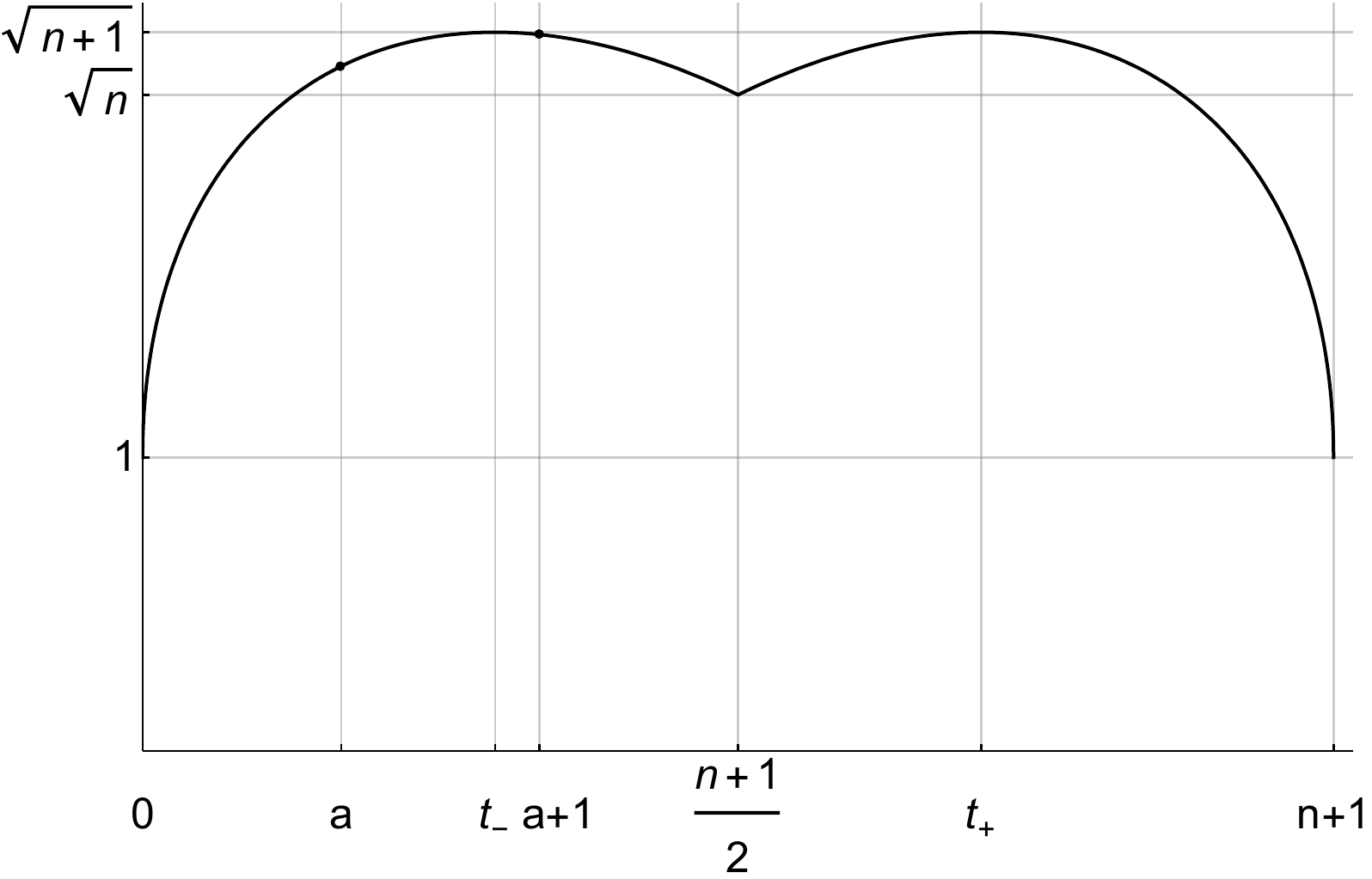}
\caption{The graph of $\psi(t)$ for $n=5$. Here $n+1=5$, $a=1$, $ t_-=
(6-\sqrt{6})/2$}
\label{fig:nev_ukl_n5}
\end{figure}

\medskip
The graph of the function $\psi(t)$ for $n=5$ is given in Fig.\ref{fig:nev_ukl_n5}. In this case $n+1=5$, $a=1$, 
and $ t_-=
(6-\sqrt{6})/2$.

Obviously, $\theta_n(B)\leq ||P||_B\leq \sqrt{n+1}$. Combining these inequalities with the inequality 
\linebreak $\theta_n(B_n)\geq c\sqrt{n}$
(see
\eqref{theta_n_B_n_gt_c_sqrt_n_1}),
 we obtain the exact order of the value of $\theta_n(B_n)$ related to dimension $n$, namely, 
$\theta_n(B_n)\asymp \sqrt{n}.$

Another corollary of Theorem 7.1 provides the following sharp values of $\theta_n(B_n)$ for $1\leq n\leq 4$:
\begin{equation}\label{theta_1234_equalities}
\theta_1(B_1)=1, \quad
\theta_2(B_2)=\frac{5}{3}, \quad
\theta_3(B_3)=2, \quad
\theta_4(B_4)=\frac{11}{5}.
\end{equation}
Of course, $\theta_1 (B_1) = 1 $, but this case also fits into the general scheme. For proving, let
us calculate the~norm of the projector
corresponding to a regular simplex inscribed into $B_n$. By  Theorem 7.1,
$$\|P\|_{B_1}=1, \quad
\|P\|_{B_2}=\frac{5}{3}, \quad
\|P\|_{B_3}=2, \quad
\|P\|_{B_4}=\frac{11}{5}.$$
Note that for each $n=1, 2, 3, 4,$ we have
$$\|P\|_{B_n}=3-\frac{4}{n+1}.$$
But $\theta_n(B_n)\geq 3-4/(n+1)$ for every positive integer $n$ so that, if  $1\leq n\leq 4$, then
\begin{equation}\label{teta_eq_3_minus_frac}
\theta_n(B_n)= 3-\frac{4}{n+1}.
\end{equation}

This is equivalent to (\ref{theta_1234_equalities}).
 If (\ref{teta_eq_3_minus_frac}) is true, then any projector having minimum norm corresponds to a regular simplex inscribed into $B_n$.
Really, thanks to \eqref{nev_ksi_P_ineq}, for any projector $P$ and the simplex $S=S_P$
$$\frac{n+1}{2}\Bigl( \|P\|_{B_n}-1\Bigr)+1\geq \xi(B_n;S)\geq \xi_n(B_n)=n,$$
i.\,e.,
$$\|P\|_{B_n}\geq3-\frac{4}{n+1}.$$
If for some projector we have here an equality, then  equalities take place  in the above chain too. In~this
case, the corresponding simplex satisfies $ \xi(B_n;S)= \xi_n(B_n)=n.$ This way $S$ is a regular
simplex inscribed in the boundary sphere (see Section 3).
Remark that the same approach provides the exact values of $\theta_n(Q_n)$ for dimensions $n=1,2,3,7.$

Let us define an integer
$k=k_n$ equal to that of the numbers $a$ or $a+1$, for which $\psi(t)$
takes a larger value.
As it shown in  \cite{nev_ukh_mais_2019_26_2}, if $k=1$, then in $B_n$ there exists an $1$-point with respect to
the simplex corresponding to $P$.
For such $ n $,
\begin{equation}\label{nev_ksi_P_eq_regular}
\xi(B_n;S)=
\frac{n+1}{2}\Bigl( \|P\|_{B_n}-1\Bigr)+1
\end{equation}
(see (\ref{nev_ksi_P_ineq})).
Since $S$ is a regular simplex inscribed into $B_n$,  then
$\xi(B_n;S)=n$  and (\ref{nev_ksi_P_eq_regular})
is equivalent to the equality $\|P\|_{B_n}=3-\frac{4}{n+1}$.
However, $k>1$ provided $n\ge 5$. This corresponds to the fact that  (\ref{teta_eq_3_minus_frac}) and  (\ref{nev_ksi_P_eq_regular}) hold only for $1\leq n\leq 4$.
Initially, this effect was discovered in the course of~computer experiments perfomed by Alexey Ukhalov.
Later the author gave an analytical solution to this problem. See \cite{nev_ukh_mais_2019_26_2} for details.

Numbers $k_n$  increase with  $n$, but not strictly  monotonously. For example,
$$k_1=k_2=k_3=k_4=1, \quad k_5=k_6=2, \quad k_7=k_8=k_9=3, \quad k_{10}=k_{11}=4,$$
$$
k_{12}=k_{13}=5, \quad k_{14}=k_{15}=6, \quad k_{50}=22, \quad k_{100}=45, \quad k_{1000}=485.
$$

Furthermore, $k_n\leq n/2$ provided $n\geq 2$.



\SECT{8. A theorem on a simplex and its minimal ellipsoid}{8}
\label{th_simp_min_ellipsoid}
\addtocontents{toc}{8. A theorem on a simplex and its minimal ellipsoid\hfill\thepage\par\VST}

\indent\par
In 1948,  F. John \cite{john_1948} proved that any convex body in ${\mathbb R}^n$ contains a unique ellipsoid of the maximum volume. Also he gave a characterization of those convex bodies for which the maximal ellipsoid
\linebreak  is the unit Euclidean ball $B_n$ (see, e.g., \cite{ball_1990}, \cite{ball_1997} for details).  John's theorem implies the analogous statement which characterizes a unique minimum volume ellipsoid containing a given convex body.
\smsk

We shall consider the minimum volume ellipsoid containing a given nondegenerate simplex.
\linebreak  For~brevity, such an ellipsoid will be called {\it the minimal ellipsoid}. Obviously, the minimal ellipsoid 
of~a~simplex is circumscribed around this simplex.
The center of the ellipsoid coincides with the center of gravity of the simplex.
The minimal ellipsoid of a simplex is a Euclidean ball if and only if this simplex is regular. This is equivalent to the well-known fact that the volume of~a~simplex contained in a ball is maximal iff this simplex is regular and inscribed into the ball
(see, e.\,g., \cite{fejes_tot_1964},
\cite{slepian_1969},
\cite {vandev_1992}).
\smsk

Let $S$ be a nondegenerate $n$-dimensional simplex in $\RN$ and let $E$ be the minimum volume ellipsoid containing $S$.
Let $m$, $1\leq m\leq n$, be a positive integer. To any $m$-point set of vertices of $S$, let us assign a point $y\in E$ defined as follows.
Let $g$ be the center of gravity of the $(m-1)$-dimensional face of $S$ which contains the chosen vertices, and let $h$ be the center of gravity of the $(n-m)$-dimensional face  which contains the rest $n+1- m$ vertices. Then $y$ is the intersection point of the straight line $(gh)$ with the boundary of the ellipsoid in direction from $g$ to $h$.

\smallskip
{\bf Theorem 8.1.} (\cite{nevskii_matzam_23}) {\it For every nondegenerate simplex $S\subset B_n$, there exists such a set of $m$ vertices for which $y\in B_n$.}

\smallskip
{\it Proof.} Let $x^{(j)}$, $j=1,...,n+1$, be the vertices of $S$. The center of gravity of the simplex and also the center of gravity of its minimal ellipsoid $E$ are  the point $c=({1}/({n+1}))\sum x^{(j)}.$ 
We let $r$ denote the~ratio of the distance between the center of gravity of a regular simplex 
and the center of gravity of~its $(m-1)$-dimensional face to the radius of the circumscribed sphere. It is easy to see that
$$
r=\frac{1}{m}\sqrt{m-\frac{m(m-1)}{n}}.
$$

Let $N$ be the number of $(m-1)$-dimensional
faces of an $n$-dimensional simplex; thus, $N={n+1\choose m}$.

Let $J$ be an $m$-element subset of the set $\{1,...,n+1\}$. The center of gravity of the $(m-1)$-dimensional face with vertices  $x^{(j)}, j\in J,$ is the point
$g_J=
(1/m) \sum_{j\in J} x^{(j)}$. Suppose $y_J$ coincides with the point $y$ for this set of vertices.
Then
$$y_J=c+\frac{1}{r}\left(c-g_J\right),$$
i.e., $y_J=(1/r)((r+1)c-g_J)$.
Summing up over all sets $J$, we have
$$\sum_{J}\|y_J\|^2=\sum_{J}\langle y_J, y_J\rangle=$$
$$=\frac{1}{r^2}\left( N\|c\|^2(r+1)^2+\sum_J\|g_J\|^2-2(r+1)\langle
\sum_J g_J,c\rangle \right)=$$
\begin{equation}\label{right_side_of_sum_y_j_norms}
=\frac{1}{r^2}\left( N\|c\|^2(r^2-1)+\sum_J\|g_J\|^2 \right).
\end{equation}
We took into  account the equality $\sum g_J=Nc.$ Further, we claim that

\begin{equation}\label{norm_g_equality}
\sum_J\|g_J\|^2 = \frac{1}{mn}{n\choose m}\sum_{j=1}^{n+1} \|x^{(j)}\|^2+
\frac{(m-1)(n+1)}{mn}{n+1\choose m} \|c\|^2.
\end{equation}
To obtain \eqref{norm_g_equality}, let us remark  the following. The value $\sum \|g_J\|^2$ contains, for every
$j$, exactly ${n\choose m-1}$ numbers $\|x^{(j)}\|^2$ taken with coefficient $1/m^2$, and
exactly ${n-1\choose m-2}$  numbers
$2\langle x^{(i)},x^{(j)}\rangle$, for each $i\ne j$, with the same coefficient $1/m^2$. (The~latter ones exist for $m\geq 2$; in case $m=1$  there are no pairwise products.)  The expression $\|c\|^2$ contains all numbers $\|x^{(j)}\|^2$,  for any $j$,
and  numbers $2\langle x^{(i)},x^{(j)}\rangle$, for~$i\ne j$, with~multiplier $1/(n+1)^2.$ The coefficient at
$\|c\|^2$ in the right-hand part of  \eqref{norm_g_equality} guarantees the~equality of~expressions with pairwise products $\langle x^{(i)},x^{(j)}\rangle$; the coefficient
at~$\sum  \|x^{(j)}\|^2$ is chosen so~that the difference of the left-hand part  of \eqref{norm_g_equality}
and the second item  in~the right-hand part is equal to the first item. Since
$$ N(r^2-1)=-\frac{(m-1)(n+1)}{mn}{n+1\choose m},$$
after replacing $\sum \|g_J\|^2$ in
\eqref{right_side_of_sum_y_j_norms} with the right-hand part of \eqref{norm_g_equality}, we notice that items with $\|c\|^2$ are~compensated. After these  transformations,
we get
\begin{equation}\label{means_equality}
\frac{1}{N}\sum_J\|y_J\|^2 =\frac{1}{n+1}\sum_{j=1}^{n+1} \|x^{(j)}\|^2.
\end{equation}
The inclusion $S\subset  B_n$ means that $\|x^{(j)}\|\leq 1.$ Therefore, the mean value of numbers $\|y_J\|^2$
is also $\leq 1.$ Consequently, for some set $J$ we have $\|y_J\|\leq 1,$
i.e., this point $y_J$ lies in $ B_n.$ Theorem 8.1 is proved.
\hfill$\Box$

\smallskip
The approach using in this proof was suggested and kindly communicated to the author by Arseniy Akopyan.

As a conjecture, Theorem 8.1 was formulated by the author in \cite{nevskii_mais_2021_28_2} .
Also in  \cite{nevskii_mais_2021_28_2},
the~statement of the theorem  for $m=1$ was proved
in a simpler way.
 In the~case $m=1$ the points $y$ have the~form
$y^{(j)}=2c- x^{(j)},$ $j=1,\ldots,n+1,$ where $x^{(j)}\in B_n$ are the vertices
of $S$ and $c=c(S)=c(E).$
We need to~show that there exists a vertex $x$ of the simplex such that $\|2c-x\|\leq 1$.
 Since simplex $S$ is~nondegenerate, for some vertex $x$ holds
 $\langle c,x-c\rangle\geq 0$. This implies
$$\|2c-x\|^2=\langle 2c-x,2c-x\rangle =
4\langle c,c-x\rangle +\|x\|^2\leq \|x\|^2\leq 1,$$
i.e., the vertex $x$ is suitable.


\SECT{9. The value of $\theta_n(B_n)$}{9}
\label{value_th_ball}
\addtocontents{toc}{9. The value of $\theta_n(B_n)$ \hfill\thepage\par\VST}

\indent\par In this section we present the main result of the work \cite{nevskii_matzam_23} devoted to calculation of the exact value of the constant $\theta_n(B_n)$.

\smallskip
{\bf Theorem 9.1.} {\it  Let $\psi$ and $a_n$ be the function and the number defined by the equalities
\eqref{psi_function_modulus} and \eqref{a_formula} respectively. Then, for every $n\in{\mathbb N}$ the following equality 
\begin{equation}\label{theta_equality}
\theta_n(B_n)=\max\{\psi(a_n),\psi(a_n+1)\}
\end{equation}
holds. Furthermore, an interpolation projector
$P:C(B_n)\to \Pi_1({\mathbb R}^n)$  satisfies the equality $\|P\|_{B_n}=\theta_n(B_n)$ if and only if its interpolation nodes coincide with the vertices
of a regular simplex inscribed into the boundary sphere.}

\smallskip
{\it Proof.}
First, let $P$ be an interpolation projector corresponding to a regular inscribed simplex  $S$, and~let $\lambda_j(x)$ be the basic Lagrange polynomials for this simplex. Let
$p_n=\|P\|_{B_n}$. Theorem 7.1 tells us that  $p_n=\max\{\psi(a_n),\psi(a_n+1)\}.$
The points $y\in B_n$ satisfying
\begin{equation}\label{p_n_equality}
\sum_{j=1}^{n+1} |\lambda_j(y)|=p_n
\end{equation}
 were found in \cite{nevskii_mais_2021_28_2}.
Namely, assume $k_n$ coincides with that of the numbers $a_n$ and $a_n+1$ which~delivers to $\psi(t)$ the maximum value. Let us fix $m=k_n$. Then  \eqref{p_n_equality} takes place for every point   $y$
which is~constructed as in Theorem 8.1 in case  $E=B_n$. The number of these points
$y$ is $N={n+1\choose m}$.

Now, let $P:C(B_n)\to \Pi_1({\mathbb R}^n)$ be an arbitrary projector with the nodes $x^{(j)}\in B_n$, $S$ be the simplex with  the vertices $x^{(j)}$, and $E$ be the minimum volume ellipsoid
contained $S$. We can also consider $P$ as a projector acting from $C(E)$. By arguments related to affine
equivalence, the norm $\|P\|_E$ coincides with the $(C(B_n)$-$C(B_n))$-norm of a projector corresponding to a regular simplex inscribed into $B_n$, i.e.,~$\|P\|_E=p_n.$
Moreover,
$$\|P\|_E=
\sum_{j=1}^{n+1} |\lambda_j(y)|$$
for $m=k_n$ and for each of the points $y$ noted at the beginning of the previous section (now with~respect to the minimal ellipsoid
$E$).

Theorem 8.1 tells us that at least one of these points, say $y^*$, belongs to $B_n$. Therefore,
$$\|P\|_{B_n}=\max_{x\in B_n} \sum_{j=1}^{n+1}|\lambda_j(x)|\geq\sum_{j=1}^{n+1}
 |\lambda_j(y^*)|=\|P\|_E=p_n.$$
 Thus, $\theta_n(B_n)=p_n=\max\{\psi(a_n),\psi(a_n+1)\}=\psi(k_n)$.

\msk
Let us show that if $\|P\|_{B_n}=\theta_n(B_n)$, then  the simplex with vertices at the interpolation nodes
is~inscribed into $B_n$ and regular. For this simplex $S$, some point  $y\in E$,
constructed for $m=k_n$, falls on~the~boundary sphere; otherwise we have
$\|P\|_{B_n}>p_n=\theta_n(B_n)$. This is due to the fact that $\sum |\lambda_j(x)|$
increases  monotonously as $x$ moves
in a straight line in direction from $c$ to $y$. 

Since $\|y\|=1,$ the mean value of $\|y_J\|$ for the rest
$y_J$ is also $\leq 1$ (see \eqref{means_equality}). Hence, another such point also lies in the boundary of the ball. Thus we obtain $\|y_J\|=1$ for all sets $J$ consisting of $m=k_n$ indices $j$. Now \eqref{means_equality} yields
$\|x^{(j)}\|=1,$ i.e., the simplex is inscribed into the ball. Since all the points
$y_J$ belong to the boundary sphere, the function $\sum |\lambda_j(x)|$  has maximum upon the ball at
$N={n+1\choose m}$ different points. Consequently, the simplex  $S$ is regular, and the proof of the theorem is  complete. \hfill$\Box$

\smallskip
For $1\leq n\leq 4$, equality \eqref{theta_equality} and characterization of minimal projectors were obtained in \cite{nev_ukh_mais_2019_26_2} and \cite{nevskii_mais_2021_28_2} by another methods.

\msk
{\bf Corollary 9.1.} {\it For every $n\in {\mathbb N}$, the following inequalities
$$
\sqrt{n}\leq \theta_n(B_n)\leq \sqrt{n+1}
$$
hold. Moreover, $\theta_n(B_n)=\sqrt{n}$ only for $n=1$ and
$\theta_n(B_n)=\sqrt{n+1}$ if and only if $\sqrt{n+1}$ is an integer.}

\smallskip
In fact, it was proved in \cite{nev_ukh_mais_2019_26_2} that the above relations are satisfied for
$p_n=\max\{\psi(a_n),\psi(a_n+1)\}$. Combining this property with the result of Theorem 9.1, we obtain Corollary 9.1.


\SECT{10. Interpolation by wider polynomial spaces}{10}
\label{int_wider_polyn_spaces}
\addtocontents{toc}{10. Interpolation by wider polynomial spaces \hfill\thepage\par\VST}

\indent\par Let us  briefly demonstrate  how some of the above methods can be applied to interpolation by wider polynomial spaces. This approach has been realized in \cite{nevskii_mais_2008_15_3}, \cite{nevskii_mais_2009_16_1}  and in details described in \cite{nevskii_monograph}.

Let $\Omega$ be a closed bounded subset of
${\mathbb R}^n$. Let $d\geq n+1,$ and let 
$\{\varphi_1(x),\ldots,\varphi_d(x)\}$ be a family of~pairwise different monomials
$x^\alpha=x_1^{\alpha_1}\ldots  x_n^{\alpha_n}.$
Here $x=(x_1,\ldots,x_n)
\in {\mathbb R}^n$,  $\alpha=(\alpha_1,\ldots,\alpha_n)$ $\in {\mathbb Z}_+^{n}$.
We suppose that
$\varphi_1(x)\equiv 1,$  $\varphi_2(x)=x_1, \
\ldots, \ \varphi_{n+1}(x)=x_n.$ 

By the $d$-dimensional space of polynomials in $n$ variables we mean the set
$\Pi=\lin(\varphi_1,\ldots,\varphi_d)$, i.e., the linear span of the family $\{\varphi_1(x),\ldots,\varphi_d(x)\}$.
Let us note the following important cases: 
$\Pi=\Pi_k\left({\mathbb R}^n\right)$ --
the space of polynomials of degree at most $k$\, $(k\in {\mathbb N})$,  and
$\Pi=\Pi_\alpha\left({\mathbb R}^n\right)$ --
the space of polynomials of degree $\leq \alpha_i$
in $x_i$ \  $(\alpha\in{\mathbb N}^n).$

A collection of points $x^{(1)},\ldots, x^{(d)} \in \Omega$ is called
{\it an admissible set of nodes } for~interpolation of~functions from $C(\Omega)$ with the use of polynomials
from $\Pi$  provided $\det({\bf A})\ne 0.$ This time
$\bf A$ is~$(d\times d)$-matrix
 $${\bf A} =
\left( \begin{array}{cccc}
1&\varphi_2\left(x^{(1)}\right)&
\ldots&\varphi_d\left(x^{(1)}\right)\\

\vdots&\vdots&\vdots&\vdots\\

1&\varphi_2\left(x^{(d)}\right)&
\ldots&\varphi_d\left(x^{(d)}\right)
\end{array}
\right).
$$
Assume that $\Omega$ contains such a set of nodes. For the interpolation projector  $P: C(\Omega)\to \Pi$ 
with these nodes, we have
\begin{equation}\label{Lagr_formula}
Pf(x)=\sum_{j=1}^d f\left(x^{(j)}\right)\lambda_j(x), \quad
\|P\|_{\Omega}=
 \max\limits_{x\in \Omega}
\sum_{j=1}^d |\lambda_j(x)|.
\end{equation}
Here $\lambda_j\in \Pi$ are the polynomials having the property
$\lambda_j\left(x^{(k)}\right)$ $=$
$\delta_j^k.$ Their coefficients with respect  to~the~basis  $\varphi_1,$ $\ldots,$ $\varphi_d$
form the~columns of  ${\bf A}^{-1}.$

Let us introduce the mapping
$T:{\mathbb R}^n\to {\mathbb R}^{d-1}$ defined in the way
$$y=T(x)=(\varphi_2(x),\ldots,\varphi_d(x))=
(x_1,\ldots,x_n,\varphi_{n+1}(x),\ldots,\varphi_d(x)).
$$

We will consider $T$ on the set $\Omega.$
The choice of the first monomials  $\varphi_j(x)$ implies the~invertibility of  $T.$
Let $y^{(j)}=T\left(x^{(j)}\right)$.
Then the numbers $\lambda_j(x)$ are the barycentric coordinates of the point
$y=T(x)$ with respect to the $(d-1)$-dimensional simplex with vertices $y^{(j)}$.

We set
$$
\theta(\Pi;\Omega)= \min\limits_{x^{(j)}\in \Omega} \|P\|_{\Omega},
$$
and
$$\xi_n(\Omega)=\min \left\{ \xi(\Omega;S): \,
S~ \mbox{is an $n$-dimensional simplex,} \,
\ver(S)\subset \Omega, \, \vo(S)\ne 0\right\}.
$$
The value $\xi(\Omega;S)$ is defined as for a convex  $\Omega$.

If $x^{(1)},\ldots,x^{(d)}\in\Omega$ is an admissible set of nodes to  interpolate
 functions from $C(\Omega)$ by~polynomials from $\Pi$, then
 $y^{(j)}=T\left(x^{(j)}\right)$ is an admissible set of nodes to interpolate
functions from  $C(T(\Omega))$ by~polynomials from
$\Pi_1\left({\mathbb R}^{d-1}\right).$   For the projector
$\overline{P}:C(T(\Omega))\to \Pi_1\left({\mathbb R}^{d-1}\right)$
with the nodes
$y^{(1)},\ldots,y^{(d)},$  we have
$\|\overline{P}\|_{T(\Omega)}=  \|P\|_\Omega. $
By this, when estimating the norm $\|P\|_{\Omega}$, it turns out to be possible to apply geometric inequalities for the norm
of projector $\overline{P}$ under linear interpolation   on the $(d-1)$ dimensional set $T(\Omega)$.
  Thus we obtain estimates of the norms $\|P\|_\Omega$ through
  the absorption coefficients $\xi(T(\Omega);S)$ or through
 the function $\chi_{d-1}^{-1}$ (modifications of \eqref{theta_n_K_ineq_thru_xi_minus_one}
 for $K=\conv(T(\Omega))$) and some others.

\smsk
For the projector $P:C(\Omega)\to\Pi$
with nodes $x^{(j)}$, we have
\begin{equation}\label{norm_P_xi_Omega_S_ineqs}
\frac{1}{2}\left(1+\frac{1}{d-1}\right)\left(\|P\|_{\Omega}-1\right)+1
\leq \xi(T(\Omega);S)\leq  \frac{d}{2}\left(\|P\|_{\Omega}-1\right)+1,
\end{equation}
where $S$ is the $(d-1)$-dimensional simplex with the vertices $T\left(x^{(j)}\right).$
Also we obtain the following inequalities:
\begin{equation}\label{xi_d_minus_1_T_Omega_S_theta_n_ineqs}
\frac{1}{2}\left(1+\frac{1}{d-1}\right)\left(\theta(\Pi;\Omega)-1\right)+1
\leq \xi_{d-1}(T(\Omega))\leq
\frac{d}{2}\left(\theta(\Pi;\Omega)-1\right)+1. \end{equation}
The right-hand inequality in \eqref{norm_P_xi_Omega_S_ineqs} becomes an equality, when
there exists a $1$-point $y=T(x)\in T(\Omega)$  with respect to the simplex
$S=\conv\left(y^{(1)},\ldots,y^{(d)}\right)$. This means that for the corresponding $x$
simultaneously
$ \|P\|_\Omega=\sum_{j=1}^d |\lambda_j(x)|$ and among the numbers $\lambda_j(x)$ there is the only one negative.

 For the case $\Pi=\Pi_1\left({\mathbb R}^n\right)$ and $\Omega=Q_n$ or $B_n$ see the previous text.
 Some other examples are given in  \cite{nevskii_monograph}.
The simplest case of  nonlinear  interpolation is quadratic interpolation on a segment.
An~analytical solution of the minimal projector problem with the indicated geometric approach is~given in \cite{nevskii_mais_2008_15_3}
(see also  \cite{nevskii_monograph} and  \cite{nev_ukh_mais_2023_30_}).
 Let us consider this case as an illustrative example.

\smsk
It is well known (see, e.\,g., \cite{pashkovskij}) that  the minimal norm of an interpolation projector is attained 
for~regular nodes and equals ${5}/{4}$. We claim that this result can be also obtained with the use of inequalities 
\eqref{norm_P_xi_Omega_S_ineqs} and \eqref{xi_d_minus_1_T_Omega_S_theta_n_ineqs}.
In addition, it turns out that there are infinitely many minimal projectors.

\begin{figure}[h]
  \centering \includegraphics[scale=0.4]{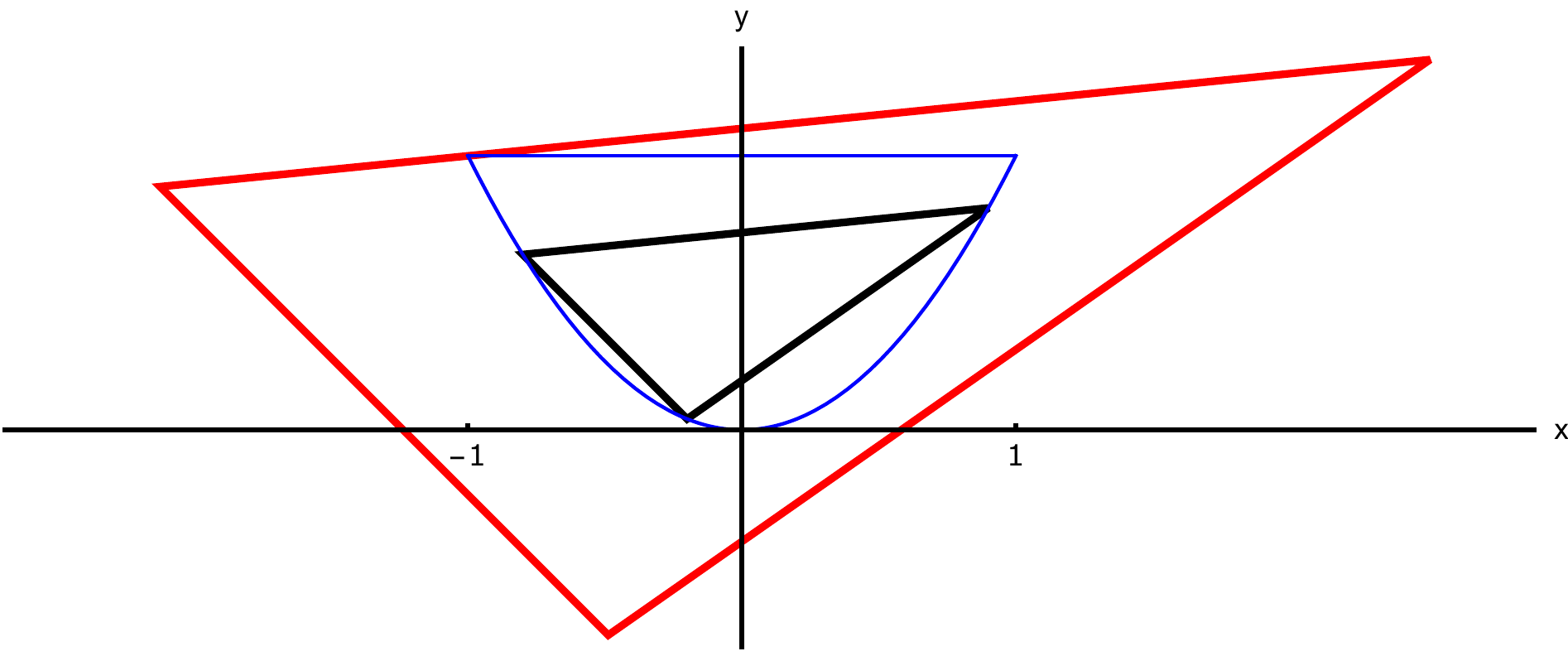}
  \caption
  {The absorption of the parabolic sector by a triangle}
  \label{fig:Triangle}
\end{figure}

\smsk
If $\Pi=\Pi_2\left({\mathbb R}^1\right)$, then $d=\dim \Pi=3,$ i.e., $k=d-1=2.$
The mapping $T$ has the form $x\longmapsto (x,x^2)$ and the set
$
T(\Omega)=T([-1,1])=\{(x,x^2)\in {\mathbb R}^2: \, -1\leq x\leq 1\}
$
is a piece of parabola. For the interpolation  nodes
$-1\leq r<s<t\leq 1$
the simplex $S$ is the triangle with  the vertices
$(r,r^2), \,(s,s^2), \,(t,t^2)$ lying on this part of parabola.
The absorption of the parabolic sector by this triangle is shown in~Fig.~\ref{fig:Triangle}\,.

The convexity of the function $\psi(x)=x^2$ implies that a $1$-point of the set  $T([-1,1])$
with respect to~$S$ exists for any nodes.
Hence, we have an equality in~the~right-hand part of \eqref{norm_P_xi_Omega_S_ineqs}:
\begin{equation}\label{norm_P_xi_eq_for_quadratic_interp}
\xi(T(\Omega);S)= \frac{3}{2}\left(\|P\|_{\Omega}-1\right)+1=\frac{3\|P\|_{\Omega}-1}{2}.
\end{equation}
Since  \eqref{norm_P_xi_eq_for_quadratic_interp} holds for an arbitrary projector
$P:C[-1,1]\to
\Pi_2({\mathbb R}^1)$, we have  the~right-hand equality also in
\eqref{xi_d_minus_1_T_Omega_S_theta_n_ineqs}:
$$\xi_2(T(\Omega))
=\frac{3\theta\left(\Pi; \Omega\right)-1}{2}.$$
So,  finding
the minimum norm of a projector
$\theta\left(\Pi; \Omega\right)$
is~equivalent to calculating
$\xi_2(T(\Omega))$.
This problem is reducible to a triangle $S$ with vertices $(-r,r^2),$ $(0,0),$ $(r,r^2),$
$0<r\leq 1.$
For~this triangle,
$$\xi(T(\Omega);S)=\max\left(\frac{11}{8},\frac{3}{r^2}-2\right), \quad
\|P\|_{\Omega}=\max\left(\frac{5}{4},\frac{2}{r^2}-1\right).$$
If
$2\sqrt{2}/3=0.9428\ldots\leq r\leq 1,$ then
these values do not depend on
$r$ and are equal correspondingly to
$11/8$ and \,$5/4$\,;
these are the minimal possible values.
Thus,
$$
\xi_2(T(\Omega))=\frac{11}{8}~~~\text{and}~~~
\theta\left(\Pi;\Omega\right)=\frac{5}{4}.
$$

The minimal is an arbitrary projector with vertices $-r,$ $0,$ $r$ where
$r\in\left[\frac{2\sqrt{2}}{3},1\right]$.
 There are no other minimal projectors.

In \cite{nev_ukh_mais_2023_30_}, we give estimates of $\theta\left(\Pi_k\left({\mathbb R}^1\right);  [-1,1]\right)$ and $\xi_k\left(T([-1,1])\right)$ for $1\leq k\leq10$. In this paper we have studied the fulfillment of equality in the right-hand side of \eqref{norm_P_xi_Omega_S_ineqs} for simplices delivering the~found upper bounds of  $\xi_k\left(T([-1,1])\right)$ and also for uniform  and Chebyshev nodes. 

This work is intended to be continued in the multidimensional case.

\bigskip
\par\noindent{\bf Acknowledgements}
\bsk

I am very thankful to Pavel Shvartsman for useful suggestions and remarks.

\bigskip


\end{document}